\documentclass[12pt]{article}
\usepackage{amsmath,amssymb,amstext,dsfont,fancyvrb,float,fontenc,graphicx,subfigure,theorem,hyperref}
\usepackage[utf8]{inputenc}
\usepackage{mathrsfs}

\usepackage{tikz-cd}
\usepackage[disable]{todonotes} 

	\usepackage{tikz}

	\usepackage[letterpaper]{geometry}
	\ifx\volno\undefined\def\volno{5}\fi
	\ifx\volyear\undefined\def\volyear{2022}\fi
	\ifx\pagno\undefined\def\pagno{24--51}\fi

	\newfont{\footsc}{cmcsc10 at 8truept}
	\newfont{\footbf}{cmbx10 at 8truept}
	\newfont{\footrm}{cmr10 at 10truept}
	
	\usepackage{fancyhdr}
	\usepackage{relsize}
	\usepackage{sectsty}
	\allsectionsfont{\larger[-1]} 
	
	\renewcommand\paragraph{\@startsection{paragraph}{4}{\z@}
		{2ex \@plus.5ex \@minus.2ex}
		{-1em}
		{\normalfont\normalsize\bfseries}}
	
	\renewcommand\subparagraph{\@startsection{subparagraph}{5}{\parindent}
		{2ex \@plus.5ex \@minus .2ex}
		{-1em}
		{\normalfont\normalsize\bfseries}}
	
	\newlength{\BiblioSpacing}
	\renewenvironment{thebibliography}[1]{
		\begin{oldthebibliography}{#1}
			\setlength{\parskip}{\BiblioSpacing}
			\setlength{\itemsep}{\BiblioSpacing}
		}
		{
		\end{oldthebibliography}
	}
	
	\usepackage[strict]{changepage}
	\def\abstractname{Abstract -}   
	\def\abstract{\begin{adjustwidth}{1cm}{1cm} \par    \footnotesize \noindent {\bf \abstractname} 
			\def\endabstract{ \end{adjustwidth} \smallskip }}
	
	{\theorembodyfont{\itshape}\newtheorem{theorem}{Theorem}[section]}
	{\theorembodyfont{\itshape}\newtheorem{proposition}[theorem]{Proposition}}
	{\theorembodyfont{\itshape}\newtheorem{definition}[theorem]{Definition}}
	{\theorembodyfont{\itshape}\newtheorem{lemma}[theorem]{Lemma}}
	{\theorembodyfont{\itshape}\newtheorem{corollary}[theorem]{Corollary}}
	
	{\theorembodyfont{\itshape}\newtheorem{question}[theorem]{Question}}
	{\theorembodyfont{\itshape}\newtheorem{alphatheorem}{Theorem}[section]}
	
	{\theorembodyfont{\itshape}}
	
	{\theorembodyfont{\rm}}
	{\theorembodyfont{\rm}}
	{\theorembodyfont{\rm}}

	\def\dedicatory{\date}
	
	\title{\Large\bf Affine Dimers from Characteristic Polygons}
	\author{\sc D. Holmes\thanks{This work was supported by the London Mathematical Society and 
	the Department of Pure Mathematics and Mathematical Statistics, University of Cambridge.}}
	\dedicatory{\normalsize\em }
	
\begin{document}
		\setcounter{page}{24}
		\maketitle
		\thispagestyle{fancy}
		
		\vskip 1.5em
		
		\begin{abstract}
			Recent work by Forsgård indicates that not every convex lattice polygon arises as the characteristic polygon of an affine dimer or, equivalently, an admissible oriented line arrangement on the torus in general position.
			We begin the classification of convex lattice polygons arising as characteristic polygons of affine dimers.
			We present several general constructions of new affine dimers from old, and an algorithm for finding affine dimers with prescribed polygon.
			
			With these tools we prove that all lattice triangles, generalised parallelograms, and polygons of genus at most two admit an affine dimer.
		\end{abstract}
		
		
		\begin{MSC}
			52C35; 52B20; 05B45 
		\end{MSC}

		
		\section{Introduction}
		\label{sec:intro}
		
		A dimer model is an embedded bipartite graph on the torus $\mathbb{T}^2$ or, depending on the application, any surface $\Sigma$.
		They were originally introduced in statistical mechanics to model molecular interactions.\todo{(B1) Explain where dimers arise in other areas of physics / maths}
		In a simplified model, the thermodynamic properties of a mixture of molecules can be calculated from a combinatorial factor that counts the number of arrangements of molecules on a square lattice.
		If all molecules are dimers rather than monomers or higher polymers, this amounts to counting the number of domino tilings of the square lattice \cite{kasteleyn1961}.\todo[disable, color=yellow]{Ising model?}
		Thinking of the square lattice as an embedded graph, there is a one-to-one correspondence between domino tilings of the lattice and perfect matchings of the graph.
		The requirement that the graph is bipartite arises when one takes into account the two possible charges of a particle.
		Finally, the torus $\mathbb{T}^2$ is a natural choice of ambient space to account for translational symmetries such as that of a crystal.
		More recent applications of dimer models can be found in algebraic and tropical geometry, as well as in string theory (e.g., \cite{addition1} and \cite{addition2}).
		
		To every dimer model on $\mathbb{T}^2$ one can associate a convex lattice polygon, called the \textit{characteristic polygon}, in at least two ways.
		The first way is as the Newton polygon of the determinant of the Kasteleyn\todo[color=green]{(A1) Kastel{\bf e}yn} operator, a generalisation of the adjacency matrix where the entries are weighted according to their meridional and longitudinal winding numbers (\cite{cimasoni2014}, Section 7).\todo{(B2) explain where the Kasteleyn matrix arises}
		A variant of this operator was used by Kasteleyn in \cite{kasteleyn1961} to calculate the number of domino tilings of a rectangular square lattice as a Pfaffian.
		The second way is as the convex hull of the values of the height function, which assigns a value in $\mathbb{Z}^2$ to every perfect matching of the dimer model. These two notions turn out to be equivalent, and it is natural to consider the inverse problem: For which convex lattice polygons does there exist a dimer model with characteristic polygon as prescribed?
		
		This question has been answered positively for all convex polygons if no further restrictions are imposed on the dimer model \cite{gulottaInverse}. Futaki–Ueda and Ueda–Yamzaki realised that, in some cases, the dimer model may be obtained from the faces of a certain hyperplane arrangement on $\mathbb{T}^2$ (\cite{futaki2010dimer}, \cite{ueda2011}, \cite{ueda2012}, as cited in \cite{forsgard2016dimer}).
		In this case, the dimer is called \textit{affine}.
		However, Forsgård exhibited a family of convex polygons which do not admit an affine dimer (\cite{forsgard2016dimer}, Section 4).
		
		The goal of this paper is to classify which convex polygons admit an affine dimer. We present partial results consisting of a list of constructions to obtain new dimers from old, an algorithm implemented in Java to verify whether a polygon admits an affine dimer, and a positive answer for all convex lattice polygons that are triangles, ``generalised parallelograms'', or have at most two interior lattice points.
		These results are summarised in Theorem \ref{thm:homologyPolygonSpace} \& \ref{thm:mainResult} at the end of this section.
		
		The results have the following application to algebraic and tropical geometry.
		\todo{(B3) explained 'no combinatorial obstruction'}
		Given a complex curve $C$ in $(\mathbb{C}^*)^2$, the \textit{coamoeba} $\mathscr{C}\subseteq \mathbb{T}^2$ is its image under the argument projection $(x,y)\mapsto (\text{arg}(x), \text{arg}(y))$ which naturally takes values on $\mathbb{T}^2$.
		The \textit{shell} of the coamoeba is a line arrangement $\mathcal{H}$ on $\mathbb{T}^2$ that is derived from the bivariate polynomial defining $C$ and satisfies $\overline{\mathscr{C}}=\mathscr{C}\cup \mathcal{H}$ (c.f. \cite{johanssonShell2013} and \cite{forsgard2016dimer}).
		Then $\mathcal{H}$ divides $\mathbb{T}^2$ into several tiles, and we say that a tile is \textit{full} if it is fully contained in $\overline{\mathscr{C}}$.
		We say that the coamoeba is represented by a dimer if we can embed a bipartite graph on $\mathbb{T}^2$ such that every vertex is contained in the interior of a full tile, every tile contains at most one vertex, and the edges correspond to shared corners between two tiles.
		If such a graph exists, it is by definition a dimer, and automatically affine since it comes from the line arrangement $\mathcal{H}$.
		
		It is natural to ask which complex curves possess coamoebas that are represented by a dimer.
		An important observation is that the Newton polygon of the defining polynomial of the curve is the same as the characteristic polygon of the affine dimer representing its coamoeba, if such a dimer exists.
		Therefore, a first obstruction is the non-existence of an affine dimer with given characteristic polygon.
		We 
		prove that this combinatorial obstruction vanishes if the genus of the curve is at most two.
		
		
		Our results also imply that all tropical curves of genus $\le 2$ can be lifted to an exact Lagrangian submanifold of $(\mathbb{C}^*)^2$, as described in \cite{jeff21}.

		For simplicity, we work with \textit{homology polygons} rather than \textit{characteristic polygons} from height functions.
		These concepts are equivalent, as outlined in the very readable source \cite{chan2016}.
		
		\subsection{Definitions}
		\label{sec:defns}
		
			The $n$\textit{-dimensional torus} $\mathbb{T}^n$ is the quotient $\mathbb{R}^n/\mathbb{Z}^n$ with quotient map $q:\mathbb{R}^n/\mathbb{Z}^n\rightarrow \mathbb{T}^n$.
			Note that $q$ is a universal cover for $\mathbb{T}^n$.
		
		\begin{definition}
			A \textit{dimer} is a bipartite multigraph $G=(V_\circ\sqcup V_\bullet, E)$ embedded on the two-dimensional torus.
			(This means we allow multiple edges between two vertices).
		\end{definition}
		
			Let $\hat H\subseteq\mathbb{R}^n$ be an affine hyperplane, i.e., there exist $a\in\mathbb{R}^n, b\in\mathbb{R}$ such that $\hat H={\{x\in\mathbb{R}^n : \langle a,x\rangle=b\}}$. We call $H:=q(\hat H)$ a \textit{hyperplane on the torus}.
			If $\text{dim}(\hat H)=1$, we call $H$ a \textit{line (on the torus)}.
		
			
			We now specialise to $n=2$.
			A \textit{closed geodesic} is a closed loop given by a line $H\subseteq\mathbb{T}^2$.
			Once we fix a choice of orientation for a closed geodesic $H$, there are unique coprime integers $a,b\in\mathbb{Z}$ such that the homology class of $H$ is $[H]=(a,b)\in H_1(\mathbb{T}^2)\cong\mathbb{Z}^2$, i.e., $(a,b)$ is the direction of $H$. We call $(a,b)\in\mathbb{Z}^2$ \textit{primitive} if $\text{\normalfont gcd}(a,b)=1$.
			One can interpret $a$ and $b$ as the winding numbers of $H$ around the two directions of $\mathbb{T}^2\cong \mathbb{T}\times \mathbb{T}$.
			Note, however, that this choice of directions is by no means intrinsic to $\mathbb{T}^2$ and we will exploit this symmetry on several occasions, using the action by automorphisms of $GL_2(\mathbb{Z})$ on $\mathbb{T}^2$.
		
\begin{figure}[h!]\centering
			\includegraphics[width=7cm]{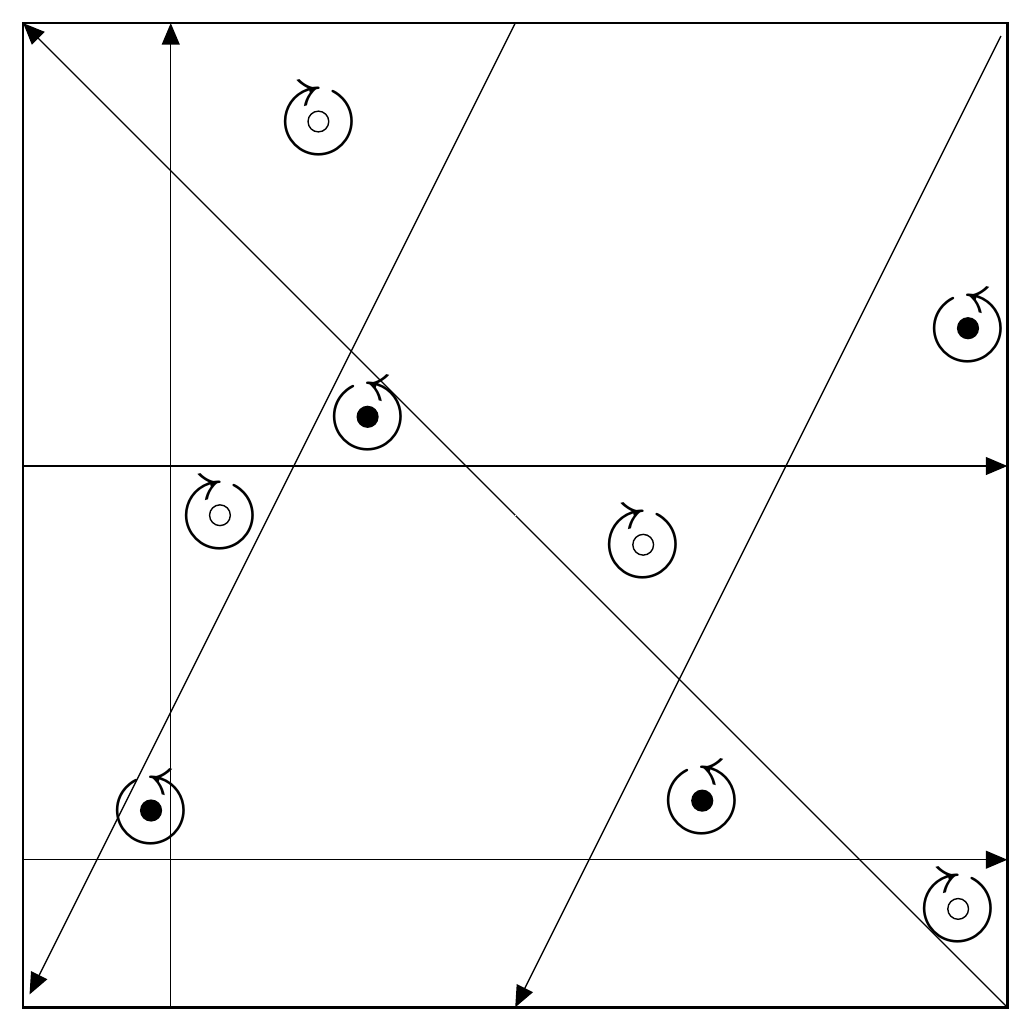}
			\hspace{1cm}
			\includegraphics[width=7cm]{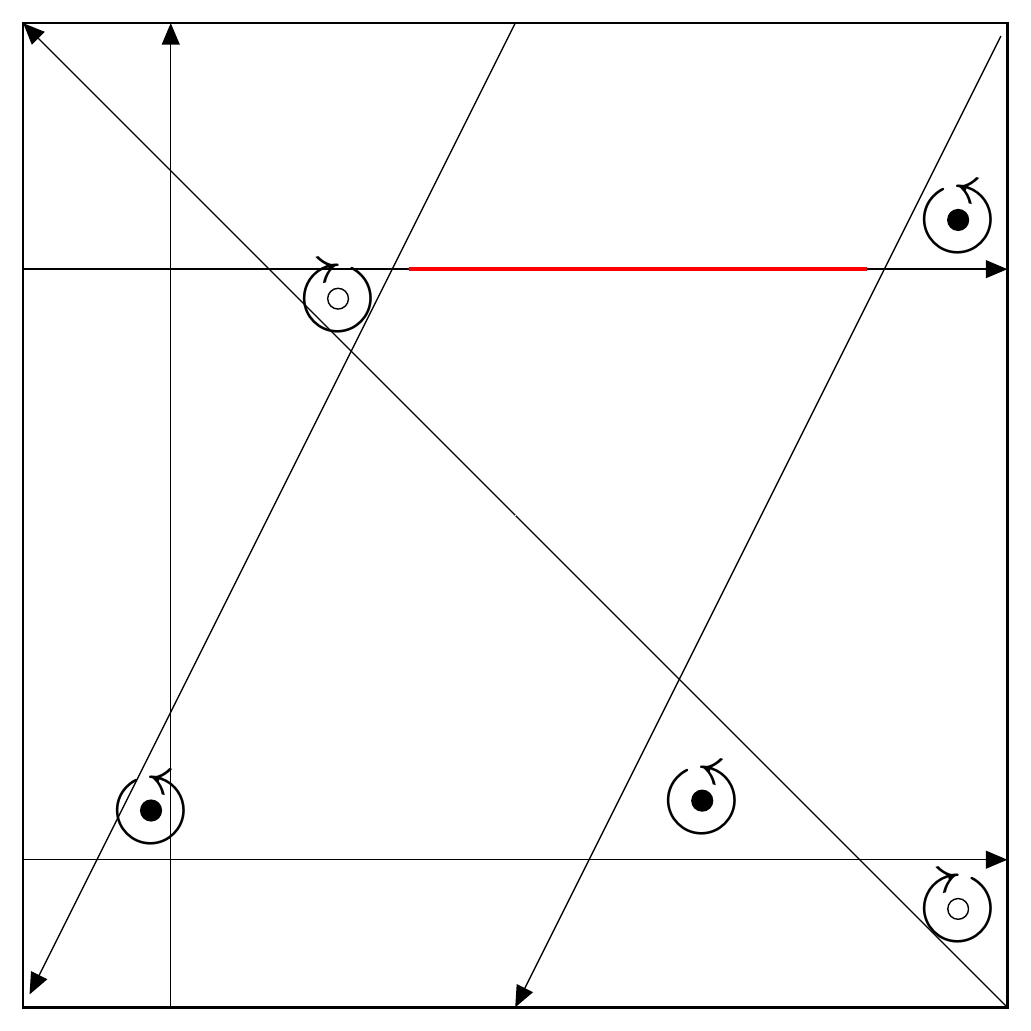}
			\caption{Examples of an admissible (left) and non-admissible (right) oriented line arrangement. Consistently oriented faces are indicated with $\circlearrowleft$ and $\circlearrowright$.
				The example to the right is not admissible because, for example, the red line segment does not bound a consistently oriented face.}
			\label{figure:admissible}
		\end{figure}

		\begin{definition}
			An \textit{oriented line arrangement} (on the 2-torus) is a finite set $\mathcal{H}$ of closed geodesics on $\mathbb{T}^2$. The line arrangement is called
			\begin{itemize}
				\item \textit{in general position} if no three lines intersect in a point, parallel lines are disjoint, and not all lines are parallel;
				\item \textit{admissible} if it is in general position and every oriented line segment is a boundary component of a face whose edges are consistently oriented, i.e., all clockwise or all counterclockwise.
				(A \textit{line segment of} $\mathcal{H}$ is a segment of a line in $\mathcal{H}$ whose endpoints are intersection points of $\mathcal{H}$ and whose interior contains no intersection points.)
			\end{itemize}
		\end{definition}
		Note that if a line arrangement is in general position then all faces are automatically homeomorphic to a disk.
		Figure \ref{figure:admissible} gives an example of\todo[color=green]{(A2) added missing {\bf of}} admissible and non-admissible line arrangements.
		Note that these only differ by a translation of the upper horizontal line, so both arrangements represent the same multiset of homology classes in $H_1(\mathbb{T}^2)$.
		
\begin{figure}[H]\centering
			\includegraphics[width=7cm]{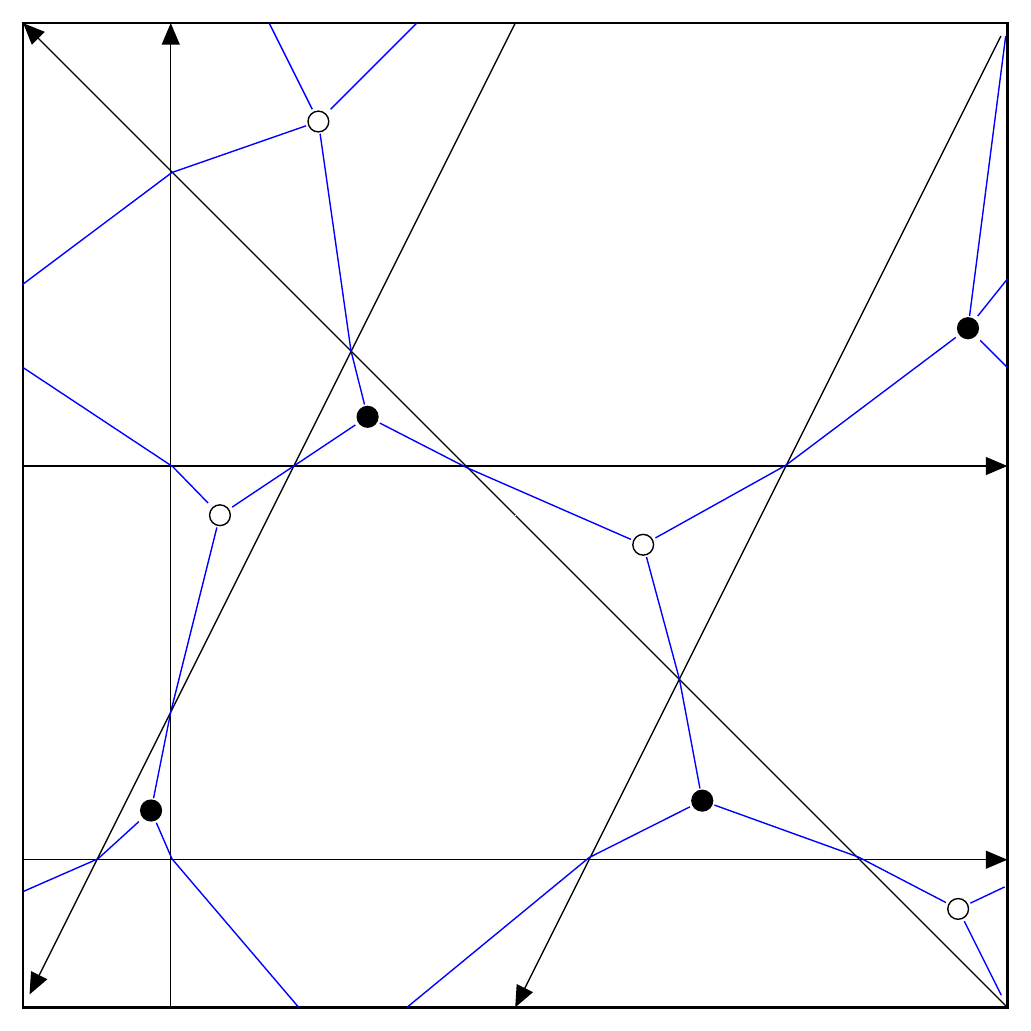}
			\caption{The affine dimer $G=(V_\circ\sqcup V_\bullet, E)$ obtained from the admissible line arrangement in Figure \ref{figure:admissible}.
				The edges of $G$ are depicted in blue.}
			\label{figure:dimer}
		\end{figure}

		We briefly elaborate on the equivalence of admissible oriented line arrangements and a certain class of dimers, called \textit{affine dimers}.
		Given an admissible oriented line arrangement, we obtain a dimer $G=(V_\circ\sqcup V_\bullet, E)$ as follows.
		Let $V_\circ$ and $V_\bullet$ be the sets of faces oriented clockwise and counterclockwise, respectively.
		For each intersection point of the line arrangement, we add an edge to $E$ connecting the two oriented faces meeting there.
		The obtained graph is bipartite.
		To embed $G$ in $\mathbb{T}^2$ we place a vertex in the interior of each consistently oriented face.
		Each edge can then be realised as a union of two line segments meeting at the shared intersection point of the two faces (see Figure \ref{figure:dimer}).
		
		The converse construction is also possible: Given a dimer $G=(V_\circ\sqcup V_\bullet, E)$ such that 
		\begin{itemize}
			\item the vertices of $G$ are faces of a line arrangement in general position;
			\item each edge of $G$ connects two faces along an intersection point of the line arrangement
			such that the connected faces are opposite each other at that intersection point;
			\item each intersection point of the line arrangement lies on exactly one edge of $G$ and each edge of $G$ contains exactly one intersection point of the line arrangement,
		\end{itemize}
		then we may declare the faces in $V_\circ$ and $V_\bullet$ to be oriented clockwise and counterclockwise, respectively.
		The above conditions determine a well-defined choice of orientation for each line. Thus, we obtain an admissible oriented line arrangement.
		
		\begin{definition}
			An \textit{affine dimer} is a dimer satisfying the three conditions above.
			Figure \ref{figure:dimer} gives an example.
		\end{definition}
		Thus, the notions of an affine dimer and an admissible oriented line arrangement are equivalent, and we will use them interchangeably.
		
		\subsection{Problem Statement}
		\label{sec:problem_statement}	
		
		Our leading question is the following:
		\begin{question}
			\label{question:main}
			For which multisets of homology classes $S=\{h_1,\dots,h_n\}\subset H_1(\mathbb{T}^2)\cong \mathbb{Z}^2$ is there an admissible oriented line arrangement $\mathcal{H}=\{H_1,\dots,H_n\}$ whose lines represent $S$, i.e., such that $[H_i]=h_i$ for $i=1,\dots,n$?
		\end{question}
		
		Figure \ref{figure:admissible} shows that a multiset of homology classes may be represented both by admissible and non-admissible oriented line arrangements.
		Moreover, Forsgård showed that there is a family of multisets of homology classes indexed by $\mathbb{N}_{\ge 5}$, for which there are no admissible oriented line arrangements representing them \cite{forsgard2016dimer}.
		Thus, the problem is non-trivial.
		
		We already saw that the homology class of a closed geodesic on $\mathbb{T}^2$ is automatically \textit{primitive}, i.e., $(a,b)\in\mathbb{Z}^2$ with $\text{\normalfont gcd}(a,b)=1$.
		There is another immediate necessary condition, which will allow us to reformulate the problem in terms of convex polygons on the integer lattice.
		
		\begin{lemma}
			\label{lemma:zeroSum}
			Let $\mathcal{H}=\{H_1,\dots,H_n\}$ be an admissible oriented line arrangement representing the homology classes $[H_i]=(a_i,b_i)\in\mathbb{Z}^2$. Then
			\[
			\sum_{i=1}^n [H_i] = 0.
			\]
		\end{lemma}
		\begin{proof}
			Each oriented line $H_i$ is subdivided into several oriented line segments whose endpoints are intersection points of the line arrangement.
			These oriented line segments represent 1-chains on $\mathbb{T}^2$, so we may write $H_i=\sum_{\text{segments } e \text{ of } H_i} e$, where $e\in C_1(\mathbb{T}^2)$ is a line segment of $H_i$.
			On the other hand, each segment belongs to exactly one consistently oriented face.
			Thus, as chains
			\[
			\sum_{i=1}^n H_i
			= \sum_{i=1}^n \left( \sum_{\text{segments } e \text{ of } H_i} e\right)
			= \sum_{\text{segments } e \text{ of }\mathcal{H}} e
			= \sum_{\text{oriented faces } F}\left(\sum_{\text{edges } e \text{ of } F} e\right).
			\]
			Passing to homology, we get $\left[\sum_{\text{edges } e \text{ of } F}e\right]=0$ for each oriented face $F$. Therefore $\sum_{i=1}^n [H_i]=0$.
		\end{proof}
		
		\begin{definition}
			A \textit{lattice polygon} is a polygon in $\mathbb{R}^2$ whose vertices all lie in $\mathbb{Z}^2$.
		\end{definition}
		
		\begin{lemma}
			There is a bijection between the finite multisets of primitive elements of $\mathbb{Z}^2$ summing to zero and the convex lattice polygons on $\mathbb{Z}^2$ up to translation.
		\end{lemma}
		\begin{proof}
			Given primitive elements $h_i=(a_i,b_i)\in\mathbb{Z}^2$, i.e., $\text{\normalfont gcd}(a_i,b_i)=1$, we may order them by their angle $\text{arg}(h_i)\in[0,2\pi)$ with the $x$-axis. We define a convex lattice polygon via the vertices $v_0=(0,0)$ and $v_i=v_{i-1}+h_i$. This is a closed polygon since $\sum_{i=1}^n h_i = 0$ and convex since we ordered the $h_i$.
			
			Conversely, given a convex lattice polygon, orient the edges counterclockwise and subdivide each edge so that it contains no integer lattice point in its interior. Viewing each edge as a vector $(a,b)\in\mathbb{Z}^2$, this corresponds exactly to $\text{\normalfont gcd}(a,b)=1$, and thus we obtain a primitive element $h_i=(a_i,b_i)\in\mathbb{Z}^2$ for each primitive edge segment of the polygon. Finally, $\sum_{i=1}^n h_i = 0$ since polygons are closed.
		\end{proof}
		
\begin{figure}[H]\centering
			\includegraphics[width=2.5cm]{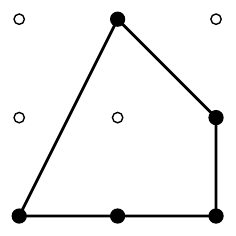}
			\caption{Homology polygon of the admissible line arrangement in Figure \ref{figure:admissible} and the affine dimer in Figure \ref{figure:dimer}. 
				The homology classes represented by the polygon are $(1,0), (1,0), (0,1),(-1,1),(-1,-2)$.}
			\label{figure:homologyPolygon}
		\end{figure}

		\begin{definition}
			Let $\mathcal{H}$ be an admissible oriented line arrangement on $\mathbb{T}^2$. The \textit{homology polygon} $P$ of $\mathcal{H}$ is the convex lattice polygon obtained from the homology classes of the lines of $\mathcal{H}$. Equivalently, we can talk about the homology polygon of an affine dimer.
			Figure \ref{figure:homologyPolygon} shows an example.
		\end{definition}
		
		It is not a priori clear that the homology polygon of an affine dimer is well-defined, as there might exist different admissible line arrangements with different homology classes that give the same affine dimer $G=(V_\circ\sqcup V_\bullet, E)$ via the construction in Section \ref{sec:defns}.
		However, as mentioned earlier, the homology polygon is equivalent to the \textit{characteristic polygon} (see \cite{chan2016}) which only depends on the data of the dimer and not the line arrangement.  Hence, the homology polygon of an affine dimer is well-defined. We only prefer to use homology polygons in the problem statement because they are slightly easier to define.
		
		Thus, we can reformulate our question:
		
		\begin{question} (Reformulation of Question \ref{question:main})
			\label{question:latticeForm}
			Which convex lattice polygons arise as the homology polygon of an admissible oriented line arrangement?
			Or which convex lattice polygons admit an affine dimer?
		\end{question}
		
		\subsubsection{Invariance under $GL_2(\mathbb{Z})$}
		Finally, we note that $GL_2(\mathbb{Z})=\left\{A\in\mathbb{Z}^{2\times 2} : \det(A)=\pm1\right\}$ acts on $\mathbb{T}^2$ by linear automorphisms which preserve admissible oriented line arrangements. Similarly, $GL_2(\mathbb{Z})$ acts on the space of convex lattice polygons through its action on $\mathbb{Z}^2$, preserving area and the number of lattice points in the interior and on the boundary. We call two lattice polygons \textit{equivalent} if they are related by an action of $GL_2(\mathbb{Z})$ and translation by a vector in $\mathbb{Z}^2$.\todo[color=green]{(A3) translation only by vectors in $\mathbb{Z}^2$.}
		Thus, whether a convex lattice polygon admits an affine dimer only depends on its equivalence class, and we arrive at our final formulation of the problem:
		
		\begin{question} (Reformulation of Question \ref{question:latticeForm})
			\label{question:latticeEquivForm}
			Which equivalence classes of convex lattice polygons arise as the homology polygon of an admissible oriented line arrangement? Or which equivalence classes admit an affine dimer?
		\end{question}
		
		\subsection{Outline of Results and Structure}
		
		Section \ref{sec:combinatorics} surveys some basic combinatorial properties of affine dimers and motivates the name \textit{genus} for the number of interior points of the homology polygon $P$, by connecting it to the genus of a punctured compact orientable surface that is homotopy equivalent to $G$.
		
		In Section \ref{sec:new_dimers_from_old} we present three constructions of affine dimers.
		The ``double everything''-construction exhibits an affine dimer for every lattice polygon consisting of pairs of antiparallel primitive side segments.
		The other two constructions (``lifting'' and ``adding an antiparallel pair'') give new dimers from old:
		\begin{alphatheorem}
			\label{thm:homologyPolygonSpace}
			Let $P$ be the homology polygon of an affine dimer.
			\begin{itemize}
				\item[(i)] If \todo[color=green]{(A4) B not b}$B\in\mathbb{Z}^{2\times 2}$ and $\det(B)\neq 0$ then $B(P)$ is also the homology polygon of an affine dimer.
				\item[(ii)] If $h\in\mathbb{Z}^2$ is a primitive side segment of $P$, then $P_h$ is also the homology polygon of an affine dimer, where $P_h$ is obtained from $P$ by adding the side segments $h,-h$.
			\end{itemize}
		\end{alphatheorem}
		\begin{proof}
			(i) Proposition \ref{prop:addParallelEdges}.
			(ii) Corollary \ref{cor:applyLinear} below.
		\end{proof}
		
		Section \ref{sec:algorithms} summarises our algorithms, including a description of the moduli space $\mathcal{M}\cong\mathbb{T}^n$ of line arrangements representing a given homology polygon $P$.
		
		Finally, Section \ref{sec:genus0and1} connects these results to finish the proof of Theorem \ref{thm:mainResult}:
		\begin{alphatheorem}\label{thm:mainResult}
			Let $P$ be a convex lattice polygon such that
			\begin{itemize}
				\item[(i)] $P$ is a triangle, or
				\item[(ii)] the primitive side segments of $P$ are pairs of antiparallel side segments, or
				\item[(iii)] the number of interior lattice points of $P$ is at most 2.
			\end{itemize}
			Then $P$ admits an affine dimer.
		\end{alphatheorem}
		\begin{proof}
			(i) Proposition \ref{prop:triangles}.
			(ii) Proposition \ref{prop:doubleEverything}.
			(iii) Propositions \ref{prop:genus0}, \ref{prop:genus1}, \ref{prop:genus2} below.
		\end{proof}
		
		\section{Basic Combinatorics of Affine Dimers}
		\label{sec:combinatorics}
		
		In this section we develop some basic combinatorics of affine dimers.
		
		\begin{definition}
			Let $G=(V_\circ\sqcup V_\bullet,E)$ be an affine dimer with corresponding admissible oriented line arrangement $\mathcal{H}=\{H_1,\dots,H_n\}$.
			Then we denote by
			\begin{itemize}
				\item $n$ ... the number of lines,
				\item $f_\circ := |V_\circ|, f_\bullet := |V_\bullet|$ ... the number of faces of the line arrangement oriented clockwise and anticlockwise, respectively,
				\item $f_\times$ ... the number of faces that are inconsistently oriented,
				\item $f=f_\circ + f_\bullet + f_\times$ ... the number of faces of the line arrangement,
				\item $v$ ... the total number of vertices of the line arrangement, i.e., intersection points of lines in $\mathcal{H}$,
				\item $e_\circ, e_\bullet$ ... the number of line segments of $\mathcal{H}$ belonging to faces in $V_\circ$ or $V_\bullet$, respectively,
				\item $e=e_\circ + e_\bullet$ ... the number of line segments of the line arrangement $\mathcal{H}$,
				\item $g$ ... the \textit{genus} of the dimer, which will be introduced in Section \ref{sec:genus}.
			\end{itemize}
		\end{definition}
		
		For example, the affine dimer in Figure \ref{figure:dimer} has $n=5, f_\circ=f_\bullet=4$, $f_\times = 5$, \\${f=v=e/2=13}$, ${e_\circ=e_\bullet=13}$, and $g=1$.
		
		\begin{proposition}[Basic counting]
			\label{prop:basicResults}
			\[
			\textnormal{(i)} \hspace{2mm}v-e+f=0\hspace{6mm}
			\textnormal{(ii)} \hspace{2mm}e_\circ = e_\bullet \hspace{6mm}
			\textnormal{(iii)} \hspace{2mm}f_\circ = f_\bullet\hspace{6mm}
			\textnormal{(iv)}\hspace{2mm} v=f=e/2
			\]
		\end{proposition}
		\begin{proof}The proofs are as follows:
			\begin{enumerate}
				\item[(i)] Immediate since $\mathbb{T}^2$ has Euler characteristic zero and a line arrangement in general position gives a CW decomposition of $\mathbb{T}^2$.
				\item[(ii)] Each of the $n$ closed geodesics consists alternately of edges counted by $e_\circ$ and $e_\bullet$.
				\item[(iii)] This follows from the existence of a perfect matching for $G$ (see Proposition \ref{lemma:matchingsExist}).
				\item[(iv)] By (i) it suffices to show $v=f$, for which we induct on the number of lines. Adding a closed geodesic in general position adds as many faces as it adds vertices.
				To verify the induction basis, assume we only have two closed geodesics which are not parallel. As the homology class of a closed geodesic is a primitive element of $\mathbb{Z}^2$, by the Euclidean algorithm we may use an action of $SL_2(\mathbb{Z})$ to assume that the geodesics have homology classes $(1,0)$ and $(c,d)$. By inspection, this configuration has $v=f=d$.
			\end{enumerate}
		\end{proof}
		
\begin{figure}[h]\centering
			\includegraphics[width=7cm]{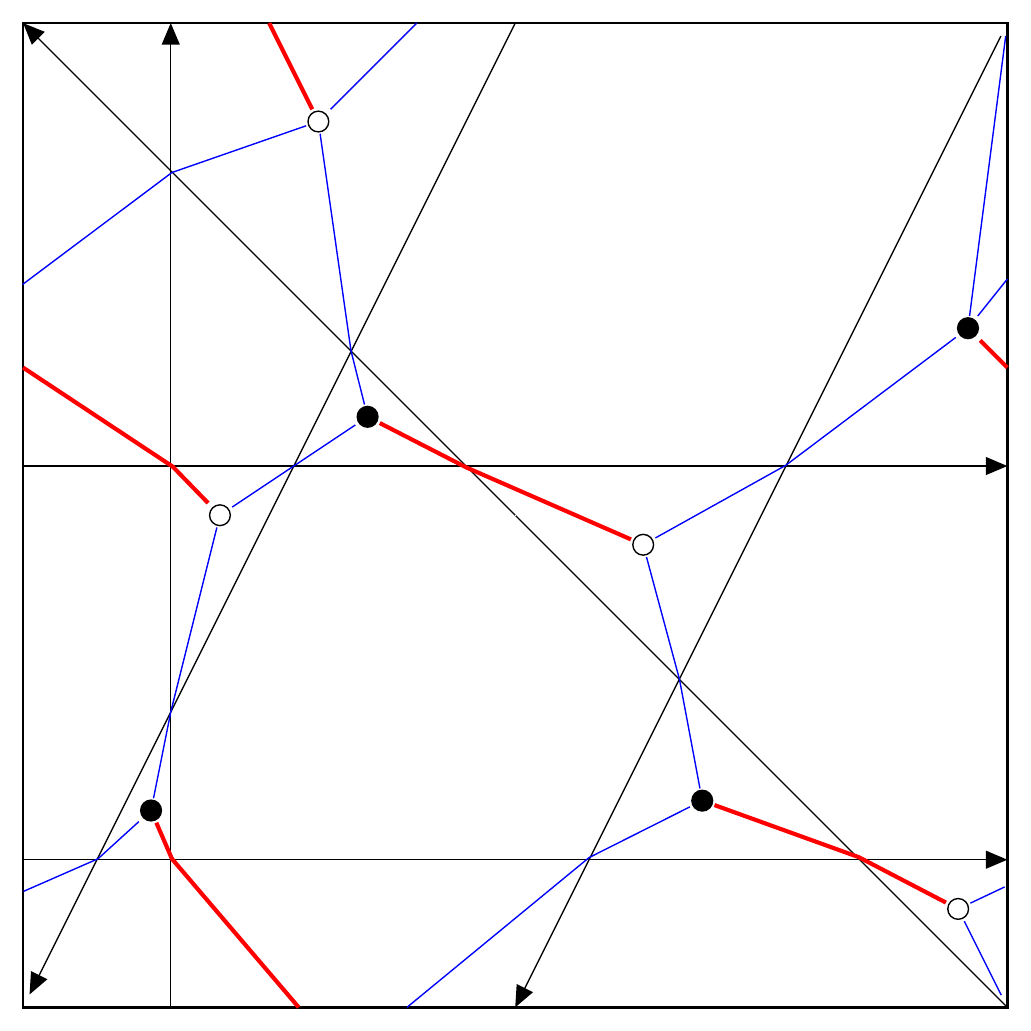}
			\hspace{1.5cm}
			\raisebox{1.1\height}{\includegraphics[width=3cm]{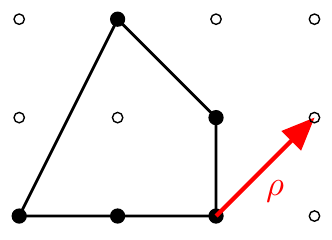}}
			\caption{Matching (red) of the affine dimer in Figure \ref{figure:dimer} corresponding to
				$\rho = (1,1)$ (left).
				As seen in the dimer's homology polygon (right), any $\rho$ with $\text{arg}(\rho)\in(0,\pi/2)$ produces the same matching.}
			\label{figure:matching}
		\end{figure}

		\begin{proposition}\label{lemma:matchingsExist}
			An affine dimer $G$ admits a perfect matching.
		\end{proposition}
		\begin{proof}
			See \cite{chan2016} for a very readable discussion of the perfect matchings of a (not necessarily affine) dimer, whose information is encoded in the \textit{characteristic polygon} via height functions. This is a special case.
			
			Let $\rho\in \mathbb{R}^2\setminus\{0\}$ be a vector that does not indicate the (signed) direction of any line in $\mathcal{H}$.
			Then every consistently oriented face has a vertex at which the directions of the two intersecting lines are immediately to the right and to the left of $\rho$.
			Now match each clockwise face in $V_\circ$ to the counterclockwise face in $V_\bullet$ adjacent to it via that vertex.
			This defines a bijection $V_\circ\rightarrow V_\bullet$ whose inverse $V_\bullet\rightarrow V_\circ$ is constructed identically using the same $\rho$. Thus, we have a matching.
			
			This construction is illustrated in Figure \ref{figure:matching}.
		\end{proof}
		
		\begin{proposition}
			\label{prop:ftimes1213}
			We have
			\[
			1/3 \le f_\times / f \le 1/2
			\]
			with $f_\times / f = 1/2$ if and only if every inconsistently oriented face is a 4-gon, and $f_\times / f = 1/3$ if and only if every consistently oriented face is a triangle.
		\end{proposition}
		\begin{proof}
			For the upper bound we count the number of corners of inconsistently oriented faces in two ways.
			On the one hand, this is $2v$ since each vertex is incident with two inconsistently oriented faces.
			On the other hand, each inconsistently oriented face has an even number of vertices, as otherwise $G$ contains an odd cycle contradicting bipartiteness. Thus,
			\[
			2v \ge 4f_\times.
			\]
			By Proposition \ref{prop:basicResults}, $v=f$, so the upper bound follows with equality condition as desired.
			
			For the lower bound we count the number of corners of consistently oriented faces in two ways. Again, this is $2v$. But every face has at least three edges, so
			\[
			2v \ge 3(f_\circ + f_\bullet) = 3(f - f_\times).
			\]
			The result follows using $f=v$ again.
		\end{proof}

		\subsection{Area of the Homology Polygon}
		
		The next result follows directly from Johansson's and Forsgård's work in \cite{forsgardJohansson2014}.\todo{(B4)=(A5)}\todo[color=green]{(A5) Index map not needed.}
		
		\begin{theorem}
			\label{thm:2areaP}
			Let $P$ be the homology polygon of an affine dimer. Then\todo[color=green]{(A6) Area not italic.}
			\[
			f_\times = 2\text{\normalfont Area}(P).
			\]
		\end{theorem}
		\begin{proofsketch}
			The key step is to show that $2\pi\text{\normalfont Area}(P)$ is the sum of \textit{inner angles} of vertices of the line arrangement, counting one per vertex  (\cite{forsgardJohansson2014}, Lemma 3.2). Here, the inner angle of two oriented intersecting lines is defined to be positive and lies between an ingoing and an outgoing ray.
			
			Thus, $2\pi\text{\normalfont Area}(P)$ is the sum of interior angles of all clockwise faces. This is exactly half the sum of exterior angles of all inconsistently oriented faces, which is $2\pi$ per face. Thus,
			\[
			2\pi\text{\normalfont Area}(P) = \frac{2\pi f_\times}{2}.
			\]
		\end{proofsketch}
		
		For example, the affine dimer in Figure \ref{figure:matching} has $f_\times = 5=2\text{\normalfont Area}(P)$.
		
		\subsection{Genus of an Affine Dimer}
		\label{sec:genus}
		
		We now describe two ways to think of an affine dimer (or equivalently of an admissible oriented line arrangement) as a two-dimensional geometric shape.
		
		\begin{definition}
			The \textit{realisation of an affine dimer} $G$ is the set $\bar{G}:=\bigcup_{F\in V_\circ\sqcup V_\bullet}\bar{F}\subseteq\mathbb{T}^2$, i.e., the union of the closed oriented faces of the admissible oriented line arrangement $\mathcal{H}$.
		\end{definition}
		
		By definition, this depends on the choice of admissible line arrangement $\mathcal{H}$ corresponding to $G$. However, many properties of $\bar{G}$ only depend on the homology polygon $P$.
		
		\begin{proposition}
			\label{prop:eulerCharFlat}
			The Euler characteristic of $\bar{G}$ is $\chi(\bar{G})=-f_\times$.
		\end{proposition}
		\begin{proof}
			\[
			\chi(\bar{G})
			= v - e + (f_\circ+f_\bullet)
			= v - e + f - f_\times = \chi(\mathbb{T}^2) - f_\times = -f_\times.
			\]
		\end{proof}
		The embedding of $G$ described in Section \ref{sec:defns} is a deformation retract of $\bar{G}$, so $\chi(G)=\chi(\bar{G})$.
		\begin{corollary}
			\label{cor:Gchar1}
			$\chi(G)=\chi(\bar{G})=-2\text{\normalfont Area}(P)$.\todo[color=green]{(A7) Area not in italic.}
		\end{corollary}
		
		We may consider $\bar{G}$ as the projection of a punctured smooth compact oriented surface $\hat{\bar{G}}$ embedded in $\mathbb{R}^3$.
		To this end we use the smooth standard embedding $\varphi: \mathbb{T}^2\hookrightarrow\mathbb{R}^3$ and consider $\varphi(\bar{G})\subset\mathbb{R}^3$.
		
\begin{figure}[H]\centering
			\includegraphics[height=3.8cm]{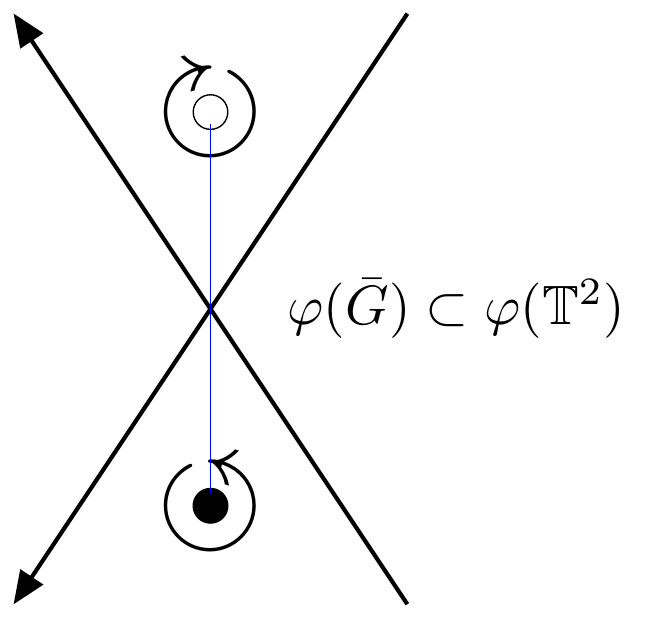}
			\hspace{1.2cm}
			\raisebox{6.5\height}{\scalebox{2}{$\longrightarrow$}}
			\hspace{1.2cm}
			\includegraphics[height=3.8cm]{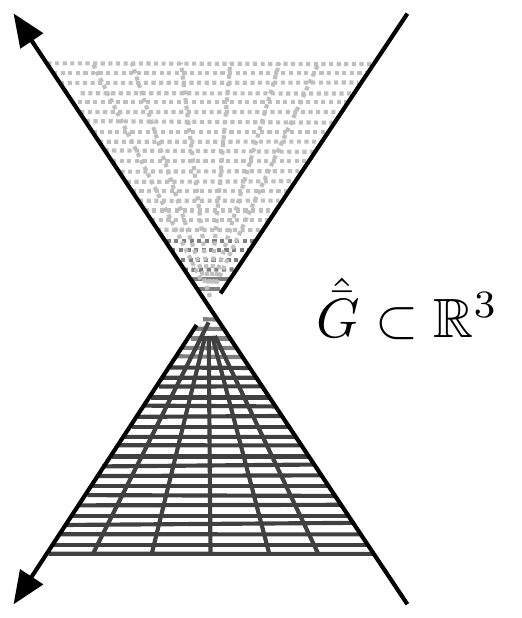}
			\caption{Twisting of $\hat{\bar{G}}$ in $\mathbb{R}^3$ near an intersection point of the boundary components of $\varphi(\bar{G})$. $\hat{\bar{G}}$ looks locally like the $180^\circ$ segment of a helicoid near such points.}
			\label{figure:twist}
		\end{figure}

		\begin{definition}(The smooth orientable surface $\hat{\bar{G}}\subset\mathbb{R}^3$)
			Away from intersection points of boundary components of $\varphi(\bar{G})$ we identify $\hat{\bar{G}}$ with $\varphi(\bar{G})$.
			Near an intersection point, we exploit the third dimension and let $\hat{\bar{G}}$ twist locally by $180^\circ$ like a helicoid as shown in Figure \ref{figure:twist}, i.e., the 
			normal vector changes smoothly from $v$ to $-v$ when traversing this neighbourhood along an edge of $G$.\todo{(B5) explain "locally like a helicoid"}
			These patches are glued together using bump functions so that we obtain a smooth compact embedded surface $\hat{\bar{G}}\subset\mathbb{R}^3$.
		\end{definition}
		
		\begin{proposition}
			\label{prop:Gorientable}
			$\hat{\bar{G}}$ is orientable.
		\end{proposition}
		\begin{proof}
			An orientation $N:\hat{\bar{G}}\rightarrow S^2$ is obtained as follows.
			By the Jordan--Brouwer\todo[color=green]{(A8) -- not -} separation theorem, $\mathbb{R}^3\setminus\varphi(\mathbb{T}^2)$ consists of a bounded and an unbounded component.
			Away from intersection points of boundary components of $\varphi(\bar{G})$, let $N(p)$ point into the unbounded component at $p$ if $p\in\varphi(\bigcup V_\bullet)$ and into the bounded component if $p\in\varphi(\bigcup V_\circ)$.
			Near an intersection point, let $N$ twist as prescribed by the local helicoid in Figure \ref{figure:twist}.
			These local definitions of $N$ glue together to form a well-defined orientation of $\hat{\bar{G}}$ because $G=(V_\circ\sqcup V_\bullet,E)$ is bipartite, so $G$ does not contain a circuit of odd length.
		\end{proof}
		
		\begin{proposition}
			\label{prop:Gchar2}
			$\bar{G}$ and $\hat{\bar{G}}$ are homotopy equivalent. Therefore,\todo[color=green]{(A7) Area not in italic.}
			\[
			\chi(\hat{\bar{G}})=\chi(\bar{G})=\chi(G)=-f_\times = -2\text{\normalfont Area}(P).
			\]
		\end{proposition}
		\begin{proof}
			It suffices to show that $\varphi(\bar{G})$ and $\hat{\bar{G}}$ are homotopy equivalent. Away from intersection points of boundary components of $\varphi(\bar{G})$, both surfaces are identical. Near an intersection point, the surfaces are equivalent by the homotopy that projects the right hand side of Figure \ref{figure:twist} onto the left hand side. These homotopies glue together compatibly and hence $\varphi(\bar{G})\simeq\hat{\bar{G}}$.
			The equalities now follow from Proposition \ref{prop:eulerCharFlat} and Corollary \ref{cor:Gchar1}.
		\end{proof}
		
		\begin{lemma}[Pick's formula]
			\label{lemma:pick}
			Let $P$ be a simple lattice polygon (i.e., $\partial P$ does not self-intersect and has exactly one connected component).\todo[color=green]{(A9) defined simple} Then\todo[color=green]{(A7) Area not in italic.}
			\[
			\text{\normalfont Area}(P) = |\mathring{P}\cap\mathbb{Z}^2| + \frac{1}{2}\left|\partial P\cap\mathbb{Z}^2\right|- 1.
			\]
		\end{lemma}
		\begin{proof}
			This is a well-known result with many different proofs available. E.g., one standard proof is via Euler's formula \cite{theBOOK}, while a more non-standard proof uses the Weierstraß $\wp$-function \cite{pickViaWeierstrass}.
		\end{proof}
		
		\begin{theorem}
			$\hat{\bar{G}}$ is homeomorphic to the compact oriented surface $\Sigma_{g,n}$ obtained by removing $n$ disjoint open discs from the compact oriented surface $\Sigma_g$ of genus $g$ without boundary.
			Moreover, the genus $g$ is the number of interior points of $P$, i.e.,
			\[
			\hat{\bar{G}}\cong\Sigma_{g,n}\hspace{5mm}\text{where}\hspace{5mm}g=\left|\mathring{P}\cap\mathbb{Z}^2\right|.
			\]
		\end{theorem}
		\begin{proof}
			The first part follows from the classification of surfaces and the fact that $\hat{\bar{G}}$ has $n$ boundary components, one for each line in $\mathcal{H}$.
			Adding a puncture to a surface (i.e., removing an open disc) decreases the Euler characteristic by one. Thus, by Proposition \ref{prop:Gchar2},
			\[
			\chi(\Sigma_g) = \chi(\hat{\bar{G}}) + n = - 2\text{\normalfont Area}(P) + n.
			\]
			But $n$ is the number of primitive side segments of $P$ which equals the number of lattice points on the boundary $\partial P$. Using $\chi(\Sigma_g)=2-2g$ we get
			\[
			g = 1 - \chi(\Sigma_g)/2 = 1 + \text{\normalfont Area}(P) - \left|\partial P\cap\mathbb{Z}^2\right|/2.
			\]
			Now the statement follows immediately from Lemma \ref{lemma:pick} (Pick's formula).
		\end{proof}
		
		This explains our definition of the genus of a dimer:
		\begin{definition}
			The \textit{genus of an affine dimer} $G$ with homology polygon $P$ is $g:=\left|\mathring{P}\cap\mathbb{Z}^2\right|$, the number of lattice points in the interior of $P$. We also call this the \textit{genus of the convex lattice polygon} $P$.
		\end{definition}
		
		This matches the relation between the genus of a tropical curve in $\mathbb{R}^2$ and the number of interior points of its Newton polygon \cite{mikhalkin2005}.\todo{(B6) reference added}
		
\begin{figure}[H]\centering
			\includegraphics[width=4cm]{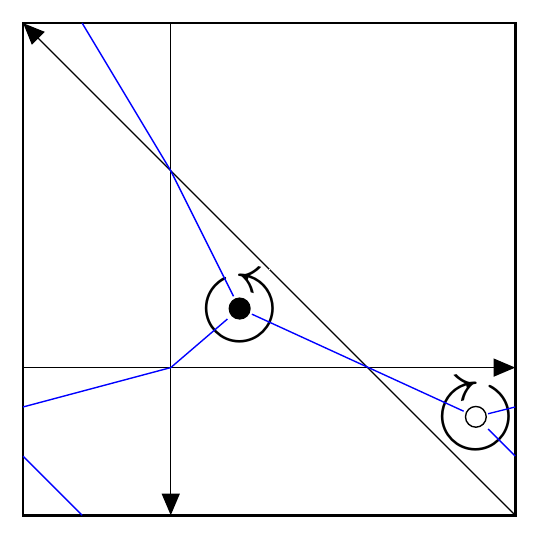}
			\hspace{1.5cm}
			\raisebox{0.7\height}{\includegraphics[width=1.7cm]{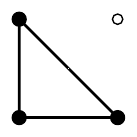}}
			\caption{An affine dimer of genus zero and its homology polygon.}
			\label{figure:genusExample}
		\end{figure}

		For example, Figure \ref{figure:genusExample} shows an affine dimer $G$ of genus zero. Since its line arrangement has $n=3$ lines, $\hat{\bar{G}}\cong \Sigma_{0,3}$ is the 3-punctured sphere also known as \textit{pair of pants}. The affine dimer in Figure \ref{figure:matching} has genus one and $\hat{\bar{G}}\cong \Sigma_{1,5}$, the 5-punctured torus.
		
		It is shown in Section \ref{sec:genus0and1} that every convex lattice polygon of genus at most 2 is the homology polygon of an affine dimer, answering Question \ref{question:latticeEquivForm} positively for these polygons.
		Moreover, by Proposition \ref{prop:triangles}, every lattice triangle admits an affine dimer. Since for every $g\in\mathbb{N}$ there exists a lattice triangle of genus $g$, there exist affine dimers of all genera.
		
		\section{Constructions of Affine Dimers}
		\label{sec:new_dimers_from_old}
		
		Next, we present three constructions of affine dimers and analyse the obtained homology polygons.
		
		\subsection{Adding parallel edges}
		
		\begin{proposition}
			\label{prop:addParallelEdges}
			Let $\mathcal{H}$ be an admissible oriented line arrangement with homology polygon $P$. Let $h\in\mathbb{Z}^2$ be a primitive side segment of $P$. Then the convex lattice polygon $P_h$ obtained by adding the antiparallel side segments $h$ and $-h$ to $P$ is the homology polygon of an admissible oriented line arrangement. Thus, if $P$ admits an affine dimer then so does $P_h$.
		\end{proposition}
		
		\begin{figure}[H]\centering
			\raisebox{0.05\height}{\includegraphics[width=6cm]{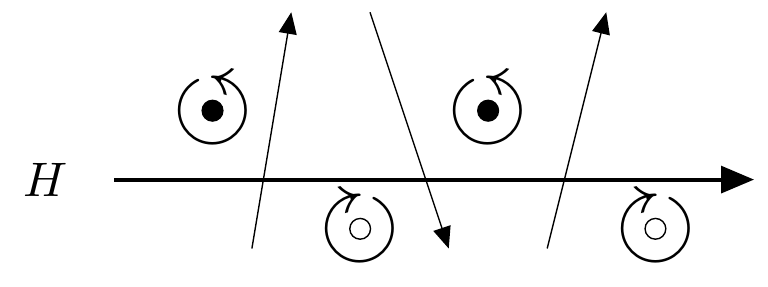}}
			\hspace{0.6cm}
			\raisebox{6\height}{\scalebox{2}{$\longrightarrow$}}
			\hspace{0.6cm}
			\raisebox{0.0\height}{\includegraphics[width=5.2cm]{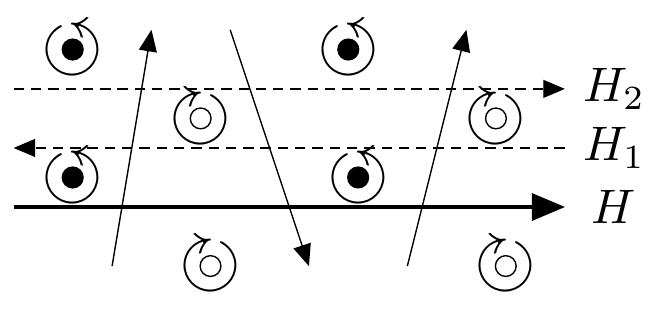}}
			\vspace{5mm}
			\newline
			\raisebox{0\height}{\includegraphics[width=3cm]{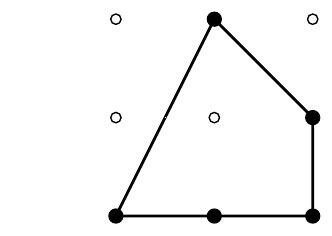}}
			\hspace{0.6cm}
			\raisebox{3\height}{\scalebox{2}{$\longrightarrow$}}
			\hspace{0.6cm}
			\raisebox{0.0\height}{\includegraphics[width=3cm]{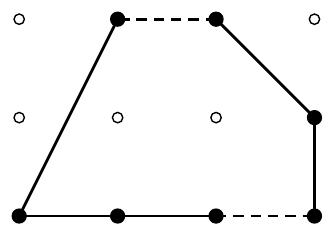}}
			\caption{Constructing $\mathcal{H}_H$ from $\mathcal{H}$ (top) and
				$P_h$ from $P$ (bottom) for $h=(1,0)$. Only the local picture near $H$ is displayed, as everything else remains unchanged.}
			\label{figure:addLine}
		\end{figure}
		
		\begin{proof}
			Let $H\in\mathcal{H}$ with $[H]=h$. We construct a new admissible oriented line arrangement $\mathcal{H}_H$ with homology polygon $P_h$ by adding two antiparallel lines $H_1,H_2$ with $[H_1]=-h$ and $[H_2]=h$.
			As depicted in Figure \ref{figure:addLine}, we place them in the order $H_2,H_1,H$ and close enough to $H$ so that no other intersection points of $\mathcal{H}$ lie between $H_2$ and $H$.
			
			If $k$ is the number of intersection points on $H$ then this construction adds $k$ consistently oriented faces to $\mathcal{H}$ locally near $H$, $k/2$ of each orientation. Away from $H$ the arrangement remains unchanged. Thus, we have obtained a new admissible arrangement $\mathcal{H}_H$ of homology polygon $P_h$, as required.
		\end{proof}

		\subsection{Double everything}
		
		All lattice polygons consisting of pairwise antiparallel primitive side segments admit an affine dimer.
		
		\begin{proposition}
			\label{prop:doubleEverything}
			Let $\Sigma=\{h_1,\dots,h_n\}\subseteq\mathbb{Z}^2$ be a multiset of primitive vectors and let $P_\Sigma$ be the convex lattice polygon consisting of the pairwise antiparallel side segments $\pm h_1, \dots,\pm h_n$.
			Then $P_\Sigma$ admits an affine dimer, i.e., $P_\Sigma$ is the homology polygon of an admissible oriented line arrangement.
			
			Additionally, the affine dimer may be taken to have $f_\times/f = 1/2$.
		\end{proposition}
		
		\begin{figure}[H]\centering
			\raisebox{0\height}{\includegraphics[width=5cm]{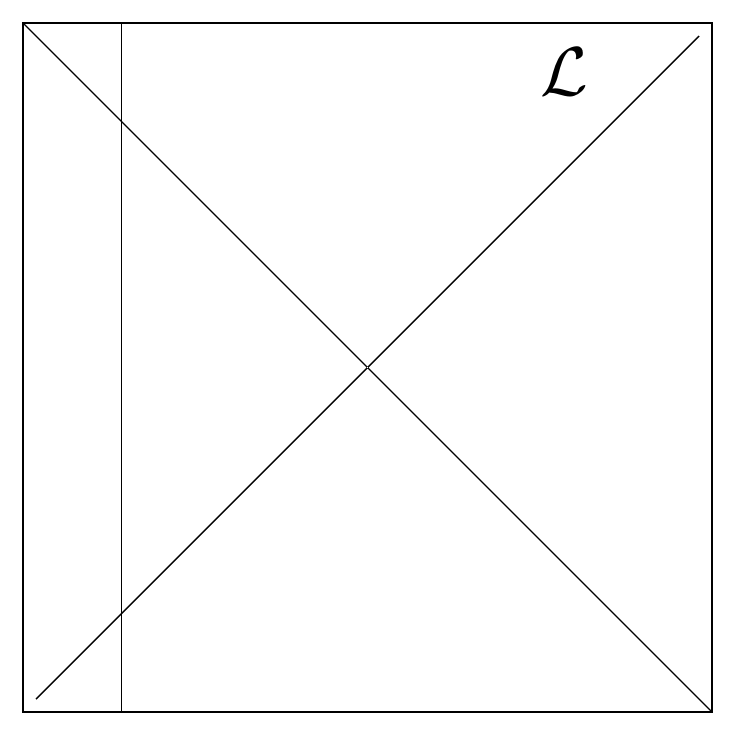}}
			\hspace{0.6cm}
			\raisebox{10\height}{\scalebox{2}{$\longrightarrow$}}
			\hspace{0.6cm}
			\raisebox{0.0\height}{\includegraphics[width=5cm]{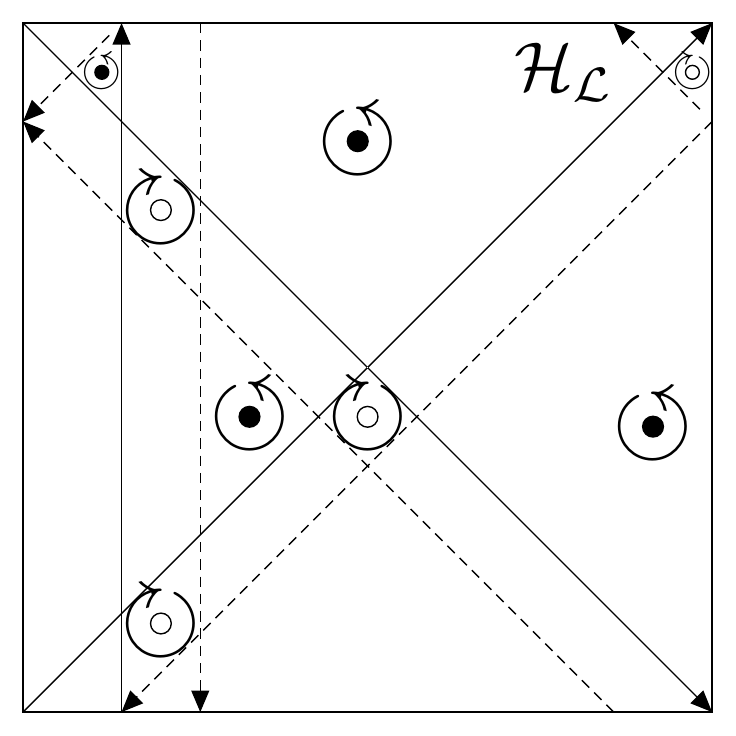}}
			\vspace{5mm}
			\newline
			\raisebox{0.75\height}{\includegraphics[height=1cm]{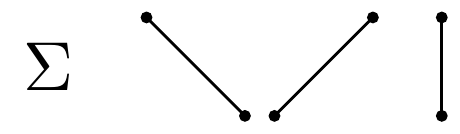}}
			\hspace{1cm}
			\raisebox{4\height}{\scalebox{2}{$\longrightarrow$}}
			\hspace{2.7cm}
			\raisebox{0\height}{\includegraphics[height=2.3cm]{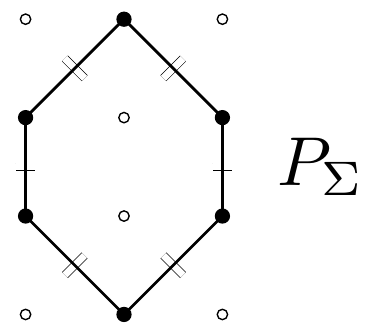}}
			\caption{Illustration of the ``double everything''-construction.}
			\label{figure:doubleEverything}
		\end{figure}
		
		\begin{proof}
			Let $\mathcal{L}$ be any unoriented line arrangement in general position representing the homology classes $\Sigma$ on $\mathbb{T}^2$ (up to sign). We construct an admissible oriented line arrangement $\mathcal{H}_\mathcal{L}$ as follows.
			First, add each line in $\mathcal{L}$ to $\mathcal{H}_\mathcal{L}$. Then, for each $H\in\mathcal{L}$, add a line $H^-$ to $\mathcal{H}_\mathcal{L}$ that is parallel and close enough to $H$, such that no lines intersect between $H$ and $H^-$ and $\mathcal{H}_\mathcal{L}$ is in general position.
			
			Let $V_\circ$ be the parallelograms corresponding to intersection points of $\mathcal{L}$ and let $V_\bullet$ be the faces corresponding to the original faces of $\mathcal{L}$.
			This gives an affine dimer $G=(V_\circ\sqcup V_\bullet, E)$ whose edges $E$ encode the face-vertex incidence relations of $\mathcal{L}$ (see Figure \ref{figure:doubleEverything}). By the discussion in Section \ref{sec:defns} there is a choice of orientation for every line in $\mathcal{H}_\mathcal{L}$ such that the faces in $V_\circ$ and $V_\bullet$ are oriented clockwise and anticlockwise, respectively. This makes  $\mathcal{H}_\mathcal{L}$ an admissible oriented line arrangement. Moreover, each pair $(H, H^-)$ is oppositely oriented, so the homology polygon of $\mathcal{H}_\mathcal{L}$ is $P_\Sigma$, as required.
			
			The inconsistently oriented faces of $\mathcal{H}_\mathcal{L}$ correspond to the line segments of $\mathcal{L}$ and are all 4-gons.
			Thus, by Proposition \ref{prop:ftimes1213}, we have $f_\times/f = 1/2$.
		\end{proof}
		
		\subsection{Lifting}
		\label{sec:lifting}
		
		In this section we use column and row vectors for elements of a vector space and its dual, respectively.
		
		Recall that $q:\mathbb{R}^2\rightarrow\mathbb{T}^2\cong\mathbb{R}^2/\mathbb{Z}^2$ is a universal cover of $\mathbb{T}^2$. Let $\mathcal{H}$ be an admissible oriented line arrangement on $\mathbb{T}^2$. The preimage $q^{-1}(\mathcal{H})$ consists of all lifts of all the lines in $\mathcal{H}$. Moreover, each fundamental parallelogram on $\mathbb{R}^2$ spanned by $\begin{pmatrix}1&0\end{pmatrix}^T, \begin{pmatrix}0&1\end{pmatrix}^T$ contains exactly one representative copy of $\mathcal{H}$.
		
\begin{figure}[H]\centering
			\raisebox{0\height}{\includegraphics[width=5cm]{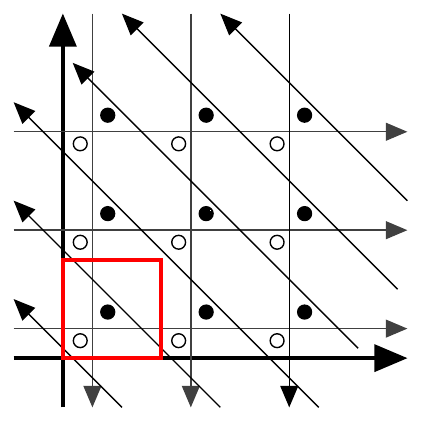}}
			\hspace{0.6cm}
			\raisebox{8\height}{\scalebox{2}{$\longrightarrow$}}
			\hspace{0.6cm}
			\raisebox{0\height}{\includegraphics[width=5cm]{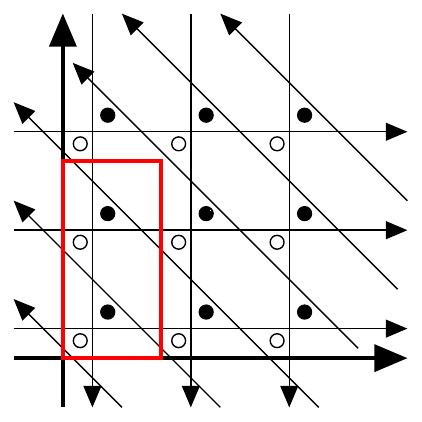}}
			\vspace{5mm}
			\\
			\raisebox{0\height}{\includegraphics[width=4cm]{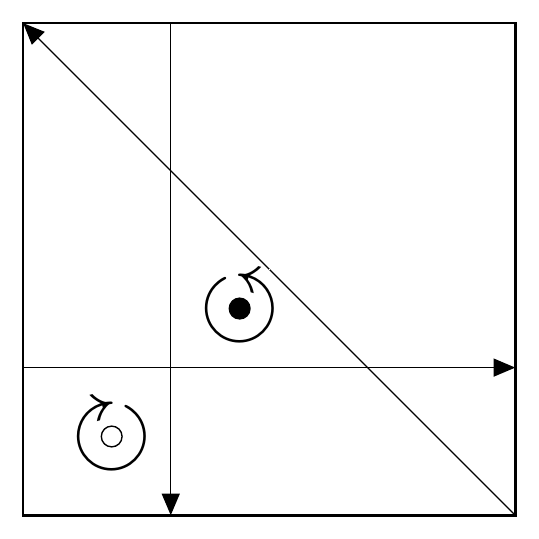}}
			\hspace{0.6cm}
			\raisebox{7\height}{\scalebox{2}{$\longrightarrow$}}
			\hspace{0.6cm}
			\raisebox{0\height}{\includegraphics[width=4cm]{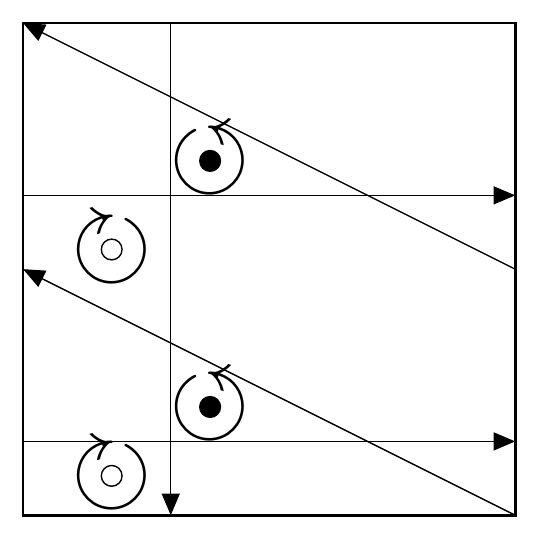}}
			\vspace{5mm}
			\\
			\raisebox{0\height}{\includegraphics[width=2.5cm]{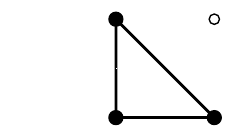}}
			\hspace{1.6cm}
			\raisebox{3\height}{\scalebox{2}{$\longrightarrow$}}
			\hspace{1.6cm}
			\raisebox{0\height}{\includegraphics[width=2.5cm]{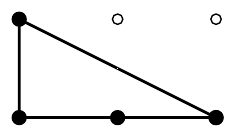}}
			\caption{The lifting construction corresponding to the lattice $\Lambda=\langle\begin{pmatrix}1&0\end{pmatrix}^T,\begin{pmatrix}0&2\end{pmatrix}^T\rangle$.
				Top: change of fundamental parallelogram. Middle: from $\mathbb{T}^2$ to $\mathbb{T}^2_\Lambda$. Bottom: the new homology polygon.}
			\label{figure:lifting1}
		\end{figure}

		However, we may define a different fundamental parallelogram spanned by two elements of $\mathbb{Z}^2$ that gives a new universal cover of a torus on which it defines a new admissible oriented line arrangement. This is equivalent to first lifting $\mathcal{H}$ to the universal cover $\mathbb{R}^2$ and then quotienting out by a general sublattice $\Lambda\le\mathbb{Z}^2$. See Figure \ref{figure:lifting1} for an example.
		
			Let $\Lambda=\langle\alpha,\beta\rangle\le\mathbb{Z}^2$ be a (non-degenerate) lattice and let $\mathbb{T}^2_\Lambda$ be the torus associated to the universal cover $q_\Lambda:\mathbb{R}^2\rightarrow\mathbb{R}^2/\Lambda:=\mathbb{T}^2_\Lambda$.
			Then $H_1(\mathbb{T}^2_\Lambda)\cong\mathbb{Z}\alpha \oplus \mathbb{Z}\beta$.
			We want to find the homology polygon of the admissible oriented line arrangement $\mathcal{H}_\Lambda:=q_\Lambda(q^{-1}(\mathcal{H}))$ on $\mathbb{T}^2_\Lambda$.
		
		Note that since $\Lambda\le\mathbb{Z}^2$, this construction gives a well-defined regular cover $q\circ q^{-1}_\Lambda:\mathbb{T}^2_\Lambda\rightarrow\mathbb{T}^2$ of degree $\left|\text{covol}(\Lambda)\right|$, the volume of any fundamental parallelogram of $\Lambda$. This cover maps $q_\Lambda(q^{-1}(\mathcal{H}))$ onto $\mathcal{H}$ so that $\mathcal{H}_\Lambda$ is a regular cover of $\mathcal{H}$.
		
		\begin{proposition}
			\label{prop:lifting}
			Let $P$ be the homology polygon of an admissible oriented line arrangement $\mathcal{H}$ on $\mathbb{T}^2$ and let $\Lambda=\langle\alpha,\beta\rangle\le\mathbb{Z}^2$ be a (non-degenerate) lattice with $\alpha=\begin{pmatrix}a\\b\end{pmatrix}$ and $\beta= \begin{pmatrix}c\\d\end{pmatrix}$.\todo{(B8) $\beta$ was missing} Let $A=\begin{pmatrix}
				a&c\\b&d
			\end{pmatrix}$ and let $P_\Lambda$ be the homology polygon of $\mathcal{H}_\Lambda$ on $\mathbb{T}^2_\Lambda$ with respect to the basis $H_1(\mathbb{T}^2_\Lambda)\cong\mathbb{Z}\alpha \oplus \mathbb{Z}\beta$.
			Then\todo[color=green]{(A10) adj not italic}
			\[P_\Lambda = \text{\normalfont adj}(A) (P)=\begin{pmatrix}
				d&-c\\-b&a
			\end{pmatrix}(P).
			\]
		\end{proposition}
		
		\begin{figure}[H]\centering
			\raisebox{0.15\height}{\includegraphics[width=3cm]{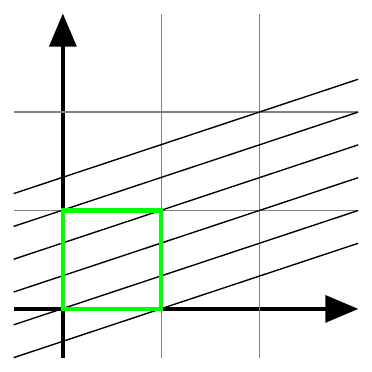}}
			\hspace{1.2cm}
			\raisebox{6\height}{\scalebox{2}{$\longrightarrow$}}
			\hspace{1.2cm}
			\raisebox{0\height}{\includegraphics[width=4cm]{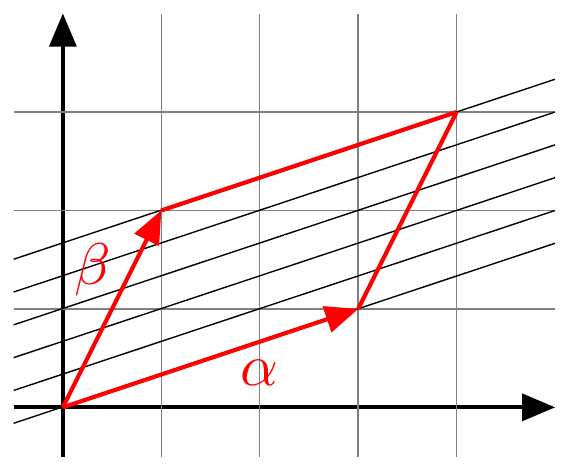}}
			\caption{Illustration of the proof of Proposition \ref{prop:lifting} with
				$\alpha=\begin{pmatrix}3&1\end{pmatrix}^T$ and $\beta=\begin{pmatrix}1&2\end{pmatrix}^T$.
				In this case $\varphi(\begin{pmatrix}3&1\end{pmatrix}^T)=\begin{pmatrix}5&0\end{pmatrix}^T$, confirming $\varphi(\alpha)=\det(A)\alpha\in H_1(\mathbb{T}^2_\Lambda)\cong\mathbb{Z}\alpha\oplus\mathbb{Z}\beta$.}
			\label{figure:lifting_proof}
		\end{figure}
		
		\begin{proof}
			Fix an orientation of $\mathbb{T}^2$.
			For two transversal loops $\gamma_1$ and $\gamma_2$ on a torus let $\iota(\gamma_1,\gamma_2)$ be their signed intersection number, which is invariant under homotopy.
			By Poincaré duality we have an isomorphism
			\begin{align*}
				i:H_1(\mathbb{T}^2)&\stackrel{\cong}{\longrightarrow} H^1(\mathbb{T}^2)\\
				[\gamma]&\longmapsto \iota([\gamma],\cdot)
			\end{align*}
			and similarly $i_\Lambda:H_1(\mathbb{T}_\Lambda^2)\stackrel{\cong}{\rightarrow} H^1(\mathbb{T}_\Lambda^2)$.
			As discussed above, the map $\pi_\Lambda:=q\circ q_\Lambda^{-1}:\mathbb{T}_\Lambda^2\rightarrow \mathbb{T}^2$ is a regular cover, restricting to a regular cover $\mathcal{H}_\Lambda$ of $\mathcal{H}$.
			Let $\varphi:=i_\Lambda^{-1}\circ \pi_\Lambda^* \circ i$.
			\begin{eqnarray}
				\label{eqn:commutative:square}
				\begin{tikzcd}
					H_1(\mathbb{T}^2) \arrow[dashed]{r}{\varphi} \arrow["i","\cong"']{d}& H_1(\mathbb{T}_\Lambda^2)\arrow["i_\Lambda","\cong"']{d}{}\\
					H^1(\mathbb{T}^2) \arrow[]{r}{\pi_\Lambda^*} &  H^1(\mathbb{T}_\Lambda^2)
				\end{tikzcd}
			\end{eqnarray}
			We shall show that $\varphi([l])=[\pi_\Lambda^{-1}(l)]$ for every $l\in Z_1(\mathbb{T}^2)$.
			By Poincaré duality this is equivalent to
			\begin{align}
				\label{eqn:poincare:iotas}
				\iota(l,\pi_\Lambda (g))=\iota(\pi_\Lambda^{-1}(l), g)
			\end{align}
			for all $l\in Z_1(\mathbb{T}^2)$ and $g\in Z_1(\mathbb{T}_\Lambda^2)$.
			Indeed, $\pi_\Lambda$ is orientation preserving and for $l:S^1\rightarrow\mathbb{T}^2$ and $g:S^1\rightarrow\mathbb{T}_\Lambda^2$ we have $l(s)=\pi_\Lambda(g(t))$ if and only if $g(t)=\hat{l}(s)$ for some lift $\hat{l}$ of $l$ along $\pi_\Lambda$, proving (\ref{eqn:poincare:iotas}).
			
			It remains to show that $\varphi$ in (\ref{eqn:commutative:square}) is given by the matrix $\text{adj}(A)$.
			Working in the basis of $H^1(\mathbb{Z}^2_\Lambda)$ dual to $H_1(\mathbb{Z}^2_\Lambda)\cong\mathbb{Z}\alpha\oplus\mathbb{Z}\beta$ we obtain the desired result:
			
			\begin{align*}
				\begin{bmatrix}
					i\begin{pmatrix}1\\0\end{pmatrix}&=&\begin{pmatrix}0&1\end{pmatrix} \\
					i\begin{pmatrix}0\\1\end{pmatrix}&=&\begin{pmatrix}-1&0\end{pmatrix}
				\end{bmatrix}
				\hspace{2mm}\Longrightarrow\hspace{2mm}
				\begin{bmatrix}
					\pi_\Lambda^*\circ i\begin{pmatrix}1\\0\end{pmatrix}&=&\begin{pmatrix}b&d\end{pmatrix} \\
					\pi_\Lambda^*\circ i\begin{pmatrix}0\\1\end{pmatrix}&=&\begin{pmatrix}-a&-c\end{pmatrix}
				\end{bmatrix}
				\hspace{2mm}\Longrightarrow\hspace{2mm}
				\begin{bmatrix}
					\varphi\begin{pmatrix}1\\0\end{pmatrix}&=&\begin{pmatrix}d\\-b\end{pmatrix} \\
					\varphi\begin{pmatrix}0\\1\end{pmatrix}&=&\begin{pmatrix}-c\\a\end{pmatrix}
				\end{bmatrix}
			\end{align*}
		\end{proof}

		The map $A\mapsto \text{adj}(A)$ on $\left\{A\in\mathbb{Z}^{2\times 2} : \det(A)\neq 0\right\}$ is surjective. Thus:
		
		\begin{corollary}
			\label{cor:applyLinear}
			If $P$ is the homology polygon of an affine dimer and $B\in\mathbb{Z}^{2\times 2}$ with $\det(B)\neq 0$, then $B(P)$ is also the homology polygon of an affine dimer.
		\end{corollary}

		We already knew this for $B\in GL_2(\mathbb{Z})$ because $GL_2(\mathbb{Z})$ acts by linear automorphisms on $\mathbb{T}^2$. This corollary is a generalisation.
		
		\section{Affine Dimer Search Algorithm}
		\label{sec:algorithms}
		
		The class of homology polygons obtained from the constructions in Section \ref{sec:new_dimers_from_old} is not too big.
		Indeed, if $P$ is a homology polygon obtained from Proposition \ref{prop:addParallelEdges} or Proposition \ref{prop:doubleEverything} then $P$ has a pair of antiparallel side segments.
		If $P$ is obtained by lifting using a matrix $B\in\mathbb{Z}^{2\times 2}$ with $\det(B)\neq 0$ as in Corollary \ref{cor:applyLinear}, and if the non-primitive side segments of $P$ are $p_1,\dots,p_m\in\mathbb{Z}^2$, then $\det(B) | \det(p_i,p_j)$ for all $i,j$.
		For this construction to deliver a new $GL_2(\mathbb{Z})$ equivalence class of convex lattice polygons, we require $\det(B)\neq \{0,\pm 1\}$.
		Thus, the integers $\det(p_i,p_j)$ all have a common prime factor, which is a rare trait for a convex lattice polygon $P$.
		
		Therefore, we developed a computer program with GUI to manipulate line arrangements on the torus and check whether a given convex lattice polygon admits an affine dimer.
		This section summarises the algorithms used.
		
		The programming was done primarily in Java, using the library \texttt{JGraphT} \cite{jgrapht} for standard graph algorithms and \texttt{polymake} \cite{polymake2020} to work with cell decompositions of $\mathbb{R}^n$.
		
		\subsection{Checking a single arrangement}
		
		Given a line arrangement $\mathcal{H}=\{H_1,\dots,H_n\}$ in general position on $\mathbb{T}^2$ with $\sum_{i=1}^n [H_i]=0$, the following algorithm determines whether it corresponds to an affine dimer.
		
\begin{figure}[h]\centering
			{\includegraphics[width=12cm]{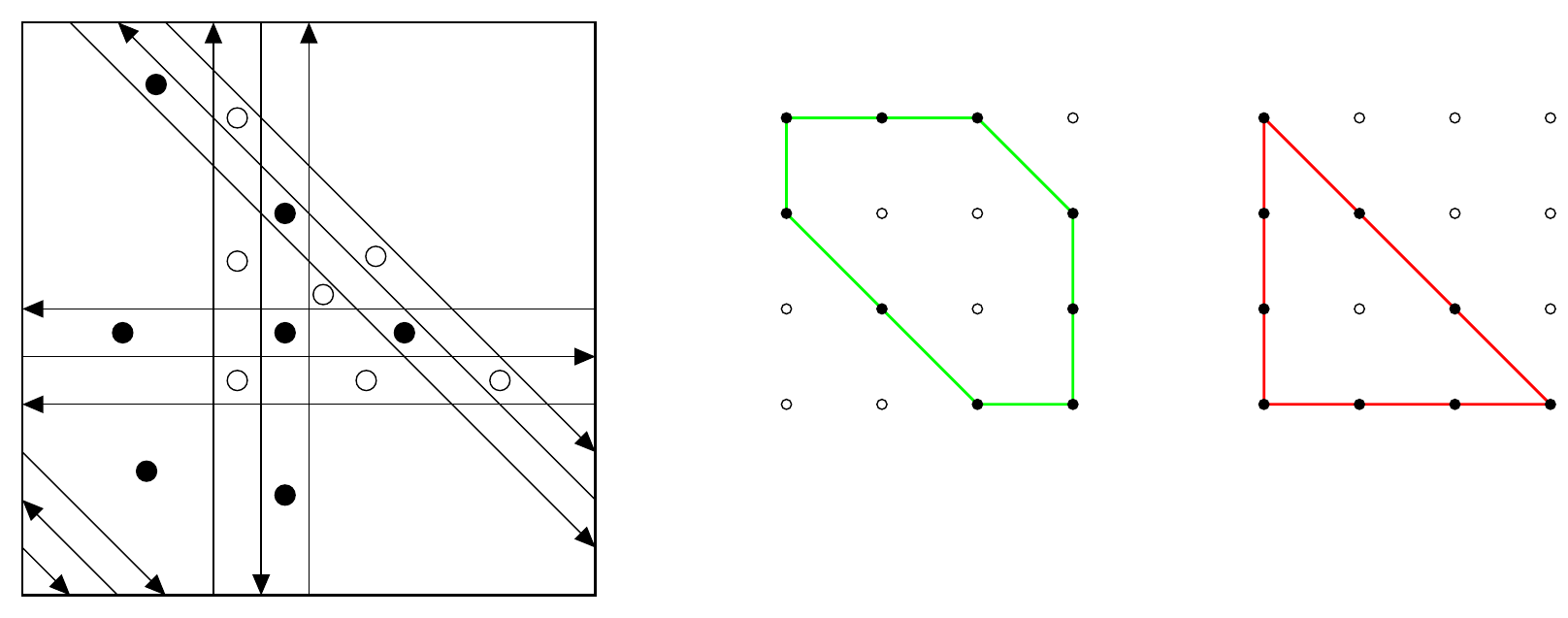}}
			\caption{The admissible arrangement to the left appears in the moduli space of the red (right) homology polygon and has $k=1$ in the above algorithm. However, it represents the green (middle) homology polygon.}
			\label{figure:pseudoDimer}
		\end{figure}

		\begin{enumerate}
			\item Calculate all intersection points. For a pair $(H_1,H_2)$ of lines, this is done by setting $[H_1]=(1,0)$ via an action of $SL_2(\mathbb{Z})$. This simplified configuration is dealt with by inspection. The number of intersection points of $H_1$ and $H_2$ is $|\det([H_1],[H_2])|$.
			
			\item For each line $H\in\mathcal{H}$, determine the order of the intersection points on $H$. Again, this is done by first setting $[H_1]=(1,0)$ via $SL_2(\mathbb{Z})$. Thus, we obtain the side segment data of $\mathcal{H}$.
			
			\item For each intersection point of each side segment, determine the next side segment at that point in clockwise and anticlockwise order. Thus, we obtain the face data of $\mathcal{H}$. Abstract this to a graph structure in which two faces are neighbours if and only if they share a vertex.
			
			\item Determine the number $k$ of bipartite connected components of the obtained graph.
			
			\item The obtained graph has exactly two connected components, which can be seen by considering intersection numbers modulo 2. Hence, there are three cases:
			\begin{itemize}
				\item If $k=2$ then the homology polygon of $\mathcal{H}$ is a parallelogram by Lemma \ref{lem:fullDimerIsParallelogram} below, which is already known to admit an affine dimer by the ``double everything''-construction of Proposition \ref{prop:doubleEverything}.
				\item If $k=1$ then there is a choice of orientation for each line in $\mathcal{H}$ making the arrangement admissible. A further check unveils whether this is compatible with the given orientations. If not, a dimer for a different homology polygon has been found.
				See Figure \ref{figure:pseudoDimer} for an example of why this is necessary.
				\item If $k=0$ then the configuration is not admissible and this cannot be fixed by re-orienting the lines.
			\end{itemize}
		\end{enumerate}
		
		This algorithm has linear time and space complexity $\mathcal{O}(f)$ by Proposition \ref{prop:basicResults} (iv), where $f$ is the number of faces of $\mathcal{H}$.
		
		\begin{lemma}
			\label{lem:fullDimerIsParallelogram}
			If $k=2$ in the above algorithm, then the homology polygon $P$ of $\mathcal{H}$ is a parallelogram.
		\end{lemma}
		\begin{figure}[H]\centering
			{\includegraphics[width=5cm]{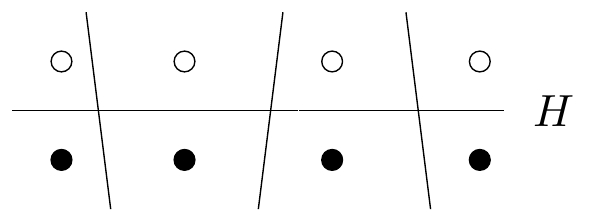}}
			\hspace{1.5cm}
			{\includegraphics[width=5cm]{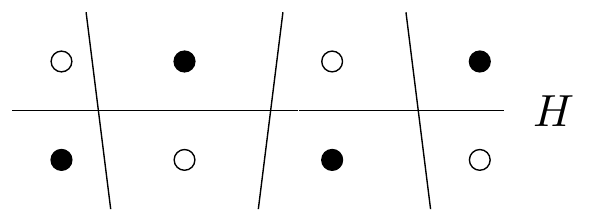}}
			\caption{The two types of lines in the case $k=2$.}
			\label{figure:fullDimer}
		\end{figure}
		\begin{proof}
			Since the obtained graph is bipartite and $k=2$, each line $H\in\mathcal{H}$ is of one of the two types depicted in Figure \ref{figure:fullDimer}. Thus, no two lines of the same type intersect, so all lines of the same type are parallel (or antiparallel). Thus, there are at most 4 homology classes and so $P$ is a parallelogram.
		\end{proof}
		
		\subsection{Moduli space of line arrangements}
		
		Checking whether a convex lattice polygon $P$ admits an affine dimer is more difficult as there are infinitely many line arrangements in general position realising $P$. However, only finitely many of them represent different combinatorial configurations. We consider two arrangements to be combinatorially the same if one of them can be obtained from the other by continuously translating some lines without ever creating a triple intersection point or two coinciding parallels. This notion is formalised by the moduli space $\mathcal{M}$ of line arrangements.
		
		\begin{lemma}
			\label{lemma:hyperplane}
			Let $\alpha\in\mathbb{Z}^n$, $c\in\mathbb{R}$, and $\hat{H}=\{x\in\mathbb{R}^n : \langle x,\alpha\rangle=c\}$ be a hyperplane. Let $q:\mathbb{R}^n\rightarrow \mathbb{T}^n$ be the quotient map and $H=q(\hat{H})$.
			Then $x+\mathbb{Z}^n\in H$ if and only if
			$\langle x,\alpha\rangle \in c + \text{\normalfont gcd}(\alpha)\mathbb{Z}$.\todo[color=green]{(A11) gcd not italic}
		\end{lemma}
		\begin{proof}
			$\mathbb{Z}$ is a Euclidean domain, so the ideal equation $(\alpha_1,\dots,\alpha_n)=\text{\normalfont gcd}(\alpha)\mathbb{Z}$ holds.
		\end{proof}
		
		Therefore, given a primitive homology class $\alpha\in\mathbb{Z}^2$, the set of lines realising this homology class is parametrised uniquely by $c\in\mathbb{R}/\mathbb{Z}\cong\mathbb{T}$, since $\text{\normalfont gcd}(\alpha)=1$.
		Therefore:
		\begin{definition}
			The \textit{moduli space} $\mathcal{M}$ of line arrangements on $\mathbb{T}^2$ consisting of $n$ lines with prescribed homology class is topologically $\mathbb{T}^n$. More precisely, if the primitive side segments of $P$ are $h_1,\dots,h_n\in\mathbb{Z}^2$ and $\alpha_i$ is the clockwise  rotation of $h_i$ by $\pi/2$ then the correspondence is
			\[
			(c_1,\dots,c_n)+\mathbb{Z}^n \in\mathbb{T}^n\cong\mathcal{M}
			\hspace{5mm}
			\longleftrightarrow
			\hspace{5mm}
			\begin{bmatrix}
				\mathcal{H}=\{H_1,\dots,H_n\},\hspace{5mm}
				[H_i] = h_i,\\
				H_i = \{x+\mathbb{Z}^2 : \langle x , \alpha_i\rangle \in c_i+ \mathbb{Z}\}
			\end{bmatrix}.
			\]
		\end{definition}
		
		Let $C\subset\mathcal{M}$ be the locus where $\mathcal{H}$ is not in general position.
		This happens either when three or more lines intersect in a point or when two (anti)parallel lines have the same parameter.
		
		For each pair $\{h_i,h_j\}$ with $i\neq j$ and $h_i \parallel h_j$, $H_i$ and $H_j$ coincide if and only if $c_i=c_j$ on $\mathbb{T}$. This happens if and only if $(c_1,\dots,c_n)$ lies on the hyperplane  $C_{i,j}: X_i-X_j=0$ on $\mathcal{M}$.
		
		For each triple $\{h_i,h_j,h_k\}$ of pairwise non-(anti)parallel homology classes, the lines $H_i,H_j,H_k$ intersect in a common point $x\in\mathbb{T}^2$ if and only if
		\[
		\langle x,\alpha_i\rangle \equiv c_i,\hspace{5mm}
		\langle x,\alpha_j\rangle \equiv c_j,\hspace{5mm}\text{and}\hspace{5mm}
		\langle x,\alpha_k\rangle \equiv c_k
		\mod \mathbb{Z}.
		\]
		This holds for some $x$ if and only if
		\[
		(c_i,c_j,c_k)\in\text{Im}(A)+\mathbb{Z}^3 \text{ for the }3\times2 \text{ integer matrix } A=\begin{pmatrix}
			\leftarrow \alpha_i \rightarrow\\
			\leftarrow \alpha_j \rightarrow \\
			\leftarrow \alpha_k \rightarrow \\
		\end{pmatrix}=:\begin{pmatrix}
			\uparrow & \uparrow \\
			e_{ijk}&f_{ijk}\\
			\downarrow&\downarrow
		\end{pmatrix},
		\]
		or equivalently $\langle (c_i,c_j,c_k), e_{ijk}\wedge f_{ijk}\rangle\in\text{\normalfont gcd}(e_{ijk}\wedge f_{ijk})\mathbb{Z}$ by Lemma \ref{lemma:hyperplane}.
		This defines a hyperplane $D_{i,j,k}\subset\mathbb{T}^n\cong\mathcal{M}$.
		Thus, the locus $C\subset\mathcal{M}$ of degenerate arrangements is\todo{(B9) only one is} given by the hyperplane arrangement
		\[
		C:=\left(\bigcup_{\substack{h_i\parallel h_j\\i\neq j}} C_{i,j}\right)\cup\left(\bigcup_{\substack{h_i,h_j,h_k\text{pairwise}\\\text{non-(anti)parallel}}} D_{i,j,k}\right)
		\subset\mathcal{M}.
		\]
		
		To determine if $P$ admits an affine dimer, one therefore needs to apply the algorithm of the previous subsection to one arrangement $(c_1,\dots,c_n)$ of each connected component of the complement $\mathcal{M}\setminus C$.
		It remains to enumerate the connected components of $\mathcal{M}\setminus C$.
		
		Note also that the coordinates of $e_{ijk}\wedge f_{ijk}$ in the definition of $D_{i,j,k}$ are of the form $\det(h_{i_1},h_{i_2})$. Thus, the structure $C\subset\mathcal{M}$ only depends on the $GL_2(\mathbb{Z})$ equivalence class the homology polygon $P$.
		
		\subsubsection{Cell decomposition approach}
		One way of enumerating the components of $\mathcal{M}\setminus C$ is as follows.
		\begin{enumerate}
			\item Lift each constituting hyperplane $H=D_{i,j,k}\text{ or }C_{i,j}$ of $C$ to several hyperplanes $\{H_1,\dots,H_l\}$ in $\mathbb{R}^n$ such that 
			$H=[0,1]^n\cap\{H_1,\dots,H_l\}\big/\mathbb{Z}^n$.
			(For $H=C_{i,j}$ we have $l=1$.
			For $H=D_{i,j,k}$ with normal vector $\nu:=e_{ijk}\wedge f_{ijk}$ we have $l\le 1+\sqrt{3}\left\Vert \nu\right\Vert / \text{\normalfont gcd}(\nu)$.)
			Thus, we obtain a (finite) hyperplane arrangement $\hat{C}\subseteq\mathbb{R}^n$ such that $q(\hat{C})=C$. 
			\item Use a cell decomposition algorithm to find the connected components of $\mathbb{R}^n\setminus\hat{C}$.
			\item Pick one point in each cell and check if the corresponding configuration is admissible.
		\end{enumerate}
		We implemented this using \texttt{polymake} \cite{polymake2020}, but the resulting algorithm was not efficient enough to deliver results for some lattice polygons $P$ of interest (such as the conjectured counterexamples of Forsgård for $k=3,4$ (\cite{forsgard2016dimer}, Section 4)).
		This is no surprise since $C$ consists of $\mathcal{O}(n^3)$ hyperplanes, each lifting to an arbitrarily large number of hyperplanes in $\hat{C}$.
		For a hyperplane arrangement in general position on $\mathbb{R}^n$ the number of components of the complement is a degree $n$ polynomial in the number of hyperplanes.
		Thus, the number of cells could be $\mathcal{O}(n^{3n})$ or even higher, depending on the number of lifts when constructing $\hat{C}$ from $C$.
		
		It might be possible to optimize this using properties of $\mathbb{T}^n$ or the fact that we are only working on $[0,1]^n$.
		
		\subsubsection{Mesh approach}
		The cell decomposition algorithm of \texttt{polymake} gives us a proper cell decomposition, whereas we only really need one point in each cell of $\mathcal{M}\setminus C$.
		There is a smallest constant $m(C)\in\mathbb{N}_{>0}$ such that every component of $\mathcal{M}\setminus C$ contains a point of $q\left(m(C)^{-1}\mathbb{Z}^n\right)$.
		Thus, we are able to finish by checking $\mathcal{O}\left(m(C)^n\right)$ configurations and without calculating any cell decompositions.
		
		However, even when $C$ just consists of two lines on $\mathbb{T}^2$, there are configurations of arbitrarily large $m(C)$.
		E.g., take lines of homology class $(1,0)$ and $(1,l)$ with $l\rightarrow\infty$.
		This configuration has $l$ faces, so $m(C)\ge l$.
		Note also that this cannot be cured by applying a smart choice of $A\in GL_2(\mathbb{Z})$ to $C$ since the number of faces is invariant.
		
		\subsubsection{Reducing dimension by two}
		In both approaches described above, we may reduce the dimension by two by restricting our attention to the subtorus $\mathbb{T}^{n-2}\subset\mathcal{M}$ with $c_1\equiv c_2\equiv 0\mod\mathbb{Z}$.
		This corresponds to translating all arrangements so that $H_1,H_2$ intersect at the bottom left corner of the fundamental parallelogram, where $h_1\not\parallel h_2$ without loss of generality.
		
		\subsection{Randomized search \& admissible volume}
		\label{sec:admissibleVolume}
		We may choose random vectors $(c_1,\dots,c_n)\in\mathcal{M}$ and check each configuration until we find an affine dimer, or stop after a certain number of trials. This approach led to the discovery of many non-trivial affine dimers of genus 1 and 2 presented in Section \ref{sec:genus0and1}.
		
		\begin{definition}
			Let $\mathcal{A}\subset\mathcal{M}\setminus C$ be the locus of admissible oriented line arrangements.
			The \textit{admissible volume} of $P$ is $\text{\normalfont vol}(\mathcal{A})$.
		\end{definition}
		
		Note that $\text{\normalfont vol}(\mathcal{A})$ only depends on the $GL_2(\mathbb{Z})$ equivalence class of the homology polygon $P$, since $C$ is $GL_2(\mathbb{Z})$ invariant. Thus, it makes sense to talk about the admissible volume of an $GL_2(\mathbb{Z})$ equivalence class of convex lattice polygons.
		
		Using this randomized approach it is possible to estimate the admissible volume of $P$.
		For some homology polygons in Section \ref{sec:genus0and1} we had $\text{\normalfont vol}(\mathcal{A})<0.01$, making it highly unlikely that we could have found their affine dimers by hand, and justifying our computational approach.
		
		Furthermore, if $\mathbb{T}^{n-2}\subset\mathcal{M}$ denotes the subtorus with $c_1\equiv c_2\equiv 0\mod\mathbb{Z}$, then $\text{\normalfont vol}_{\mathbb{T}^n}(\mathcal{A})=\text{\normalfont vol}_{\mathbb{T}^{n-2}}(\mathcal{A}\cap\mathbb{T}^{n-2})$.
		This is because the $(n-2)$-dimensional fibers are isomorphic via global translation of the arrangements, and translation is an isometry on flat tori.
		Thus, we may speed up the estimation of $\text{\normalfont vol}(\mathcal{A})$ (and thus our search for affine dimers) by reducing the dimension by two.
		
		Finding bounds for $\text{\normalfont vol}(\mathcal{A})$ in terms of $P$ might allow us to answer Question \ref{question:latticeEquivForm} for bigger classes of polygons.
		For example, we have the following result for parallelograms, which might be generalised in future.
		
		\begin{proposition}
			\label{prop:parallelogramAdmissibleVol}
			Let $P$ be an $a\times b$ lattice parallelogram with $n=2a+2b$ primitive side segments.
			Then\todo[color=green]{(A12) vol not italic}
			\[
			\text{\normalfont vol}(\mathcal{A})=4
			\begin{pmatrix}
				2a\\a
			\end{pmatrix}^{-1}
			\begin{pmatrix}
				2b\\b
			\end{pmatrix}^{-1}.
			\]
		\end{proposition}
		\begin{proof}
			We say that $P$ has $2a$ horizontal and $2b$ vertical side segments, $a$ (respectively $b$) of each orientation.
			An arrangement $(c_1,\dots,c_n)\in\mathcal{M}\setminus C$ is admissible if and only if both the vertical lines and the horizontal lines have alternating orientations.
			There are $\begin{pmatrix}2a \\ a\end{pmatrix}$ orders of orientations for the horizontal lines, exactly two of which are alternating (and similarly for the vertical lines).\todo{(B10) lines not liens}
			The result follows since all orderings of horizontal (respectively vertical) lines are equally likely and from independence of vertical and horizontal lines.
		\end{proof}
		
		\section{Triangles and Genus $\boldmath\le 2$}
		\label{sec:genus0and1}
		
		We now apply the constructions of Section \ref{sec:new_dimers_from_old} and the algorithms of Section \ref{sec:algorithms} to exhibit various families of convex lattice polygons that admit an affine dimer.
		We also record some estimates of their admissible volumes (see Section \ref{sec:admissibleVolume}), indicating how hard it would be to find their affine dimers by hand without our constructions and algorithms.
		
		\subsection{Triangles}
		
		We begin with an application of the lifting construction in Section \ref{sec:lifting}.
		
		\begin{proposition}
			\label{prop:triangles}
			Let $P$ be a lattice triangle (possibly with more than three primitive side segments).
			Then $P$ admits an affine dimer.
			Moreover, it admits an affine dimer that is lifted from the one in Figure \ref{figure:genusExample}.
		\end{proposition}
		\begin{proof}
			Let the triangle $P$ be spanned by the vectors $(a,b), (c,d)\in\mathbb{Z}^2$, not necessarily primitive as we allow more than three primitive side segments.
			Then an admissible oriented line configuration with homology polygon $P$ is obtained by applying the matrix $B=\begin{pmatrix}
				a&c\\b&d
			\end{pmatrix}$ to the triangle spanned by $(1,0)$ and $(0,1)$ using Corollary \ref{cor:applyLinear}.
		\end{proof}
		
		\subsection{Genus Zero}
		
		\begin{figure}[H]\centering
			\raisebox{0\height}{\includegraphics[width=10cm]{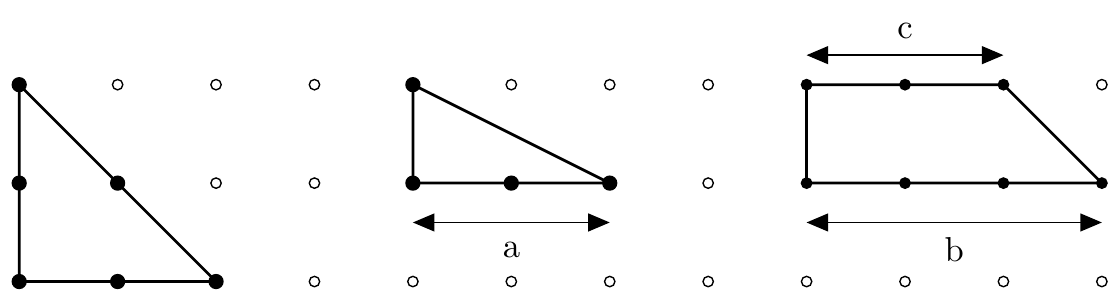}}
			\caption{The three families of equivalence classes of convex lattice polygons with no interior lattice points, where $a,b,c$ are positive integers \cite{genus0polygons}.}
			\label{figure:genus0classes}
		\end{figure}
		
		\begin{proposition}
			\label{prop:genus0}
			Let $P$ be a convex lattice polygon with no interior lattice points.
			Then $P$ admits an affine dimer.
		\end{proposition}
		\begin{proof}
			By \cite{genus0polygons} the equivalence classes of convex lattice polygons with no interior lattice points are those displayed in Figure \ref{figure:genus0classes}.
			By Proposition \ref{prop:triangles}, the triangles all admit an affine dimer.
			
			For the trapezoid given by $b,c\in\mathbb{Z}_{>0}$, there are two cases.
			If $b=c$ then $P$ consists of pairwise antiparallel edges, so admits an affine dimer by Proposition \ref{prop:doubleEverything}.
			If $b\neq c$ and without loss of generality $b>c$ then $P$ is obtained from the triangle in Figure \ref{figure:genus0classes} with $a=c$ by adding $b-c$ pairs of antiparallel edges parallel to $(1,0)$. Thus, $P$ admits an affine dimer by Proposition \ref{prop:triangles} and Proposition \ref{prop:addParallelEdges}.
		\end{proof}
		
		Even for the genus 0 triangles the admissible volume decays quickly:
		\begin{proposition}
			Let $P$ be the $a\times 1$ triangle in Figure \ref{figure:genus0classes}. The admissible volume of $\mathcal{A}\subset\mathcal{M}\setminus C$ is
			\[
			\text{\normalfont vol}(\mathcal{A})=a! / a^a.
			\]
		\end{proposition} 
		\begin{proof}
			Let $H_1,\dots,H_a$ be the lines of homology class $(1,0)$ and let $V$ and $S$ be the lines of homology class $(0,-1)$ and $(-a,1)$, respectively.
			Then $S$ subdivides $V$ into $a$ segments $s_1,\dots,s_a$ of equal length.
			Since the arrangement is in general position, there is a function $f:\{1,\dots,a\}\rightarrow\{1,\dots,a\}$ such that $H_i$ intersects $V$ in the segment $s_{f(i)}$.
			By inspection, the arrangement is admissible if and only if $f$ is injective.
			Since every $f$ is equally likely for $(c_1,\dots,c_{a+2})\in\mathcal{M}\setminus C$ chosen uniformly at random, we have
			\[
			\text{\normalfont vol}(\mathcal{A})=\mathbb{P}(f \text{ is injective}) = a! / a^a.
			\]
		\end{proof}
		
		\subsection{Genus One}
		
		\begin{figure}[H]\centering
			\includegraphics[width=14cm]{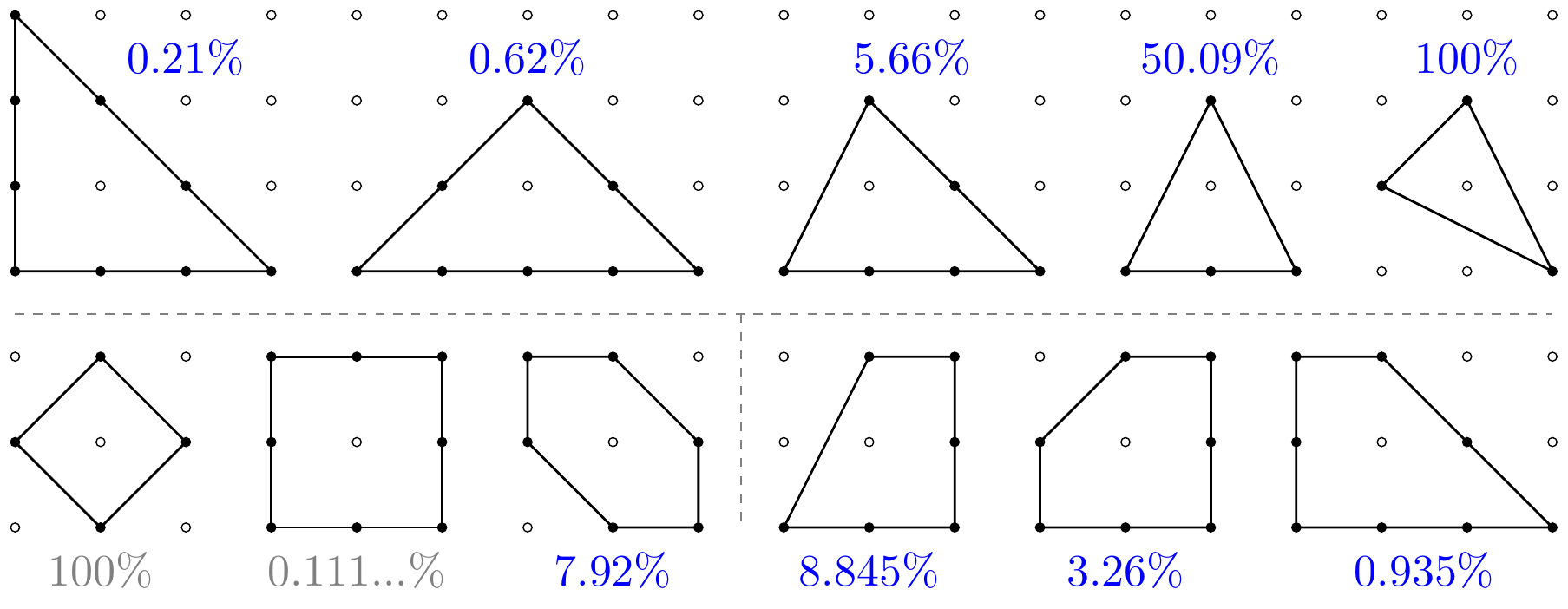}
			\caption{Equivalence class representatives of convex lattice polygons with one interior point that are
				triangles (top), consist of pairwise antiparallel side segments (bottom left), or have a pair of antiparallel side segments which are parallel to at least one other side segment (bottom right).
				The blue numbers are estimates ($\ge2\cdot 10^4$ trials) of the admissible volume of the polygon, indicating how hard it would be to find affine dimers by hand.
				The two bottom left numbers are exact by Proposition \ref{prop:parallelogramAdmissibleVol}.}
			\label{figure:genus1easyClasses}
		\end{figure}
		
		\begin{proposition}
			\label{prop:genus1}
			Let $P$ be a convex lattice polygon with exactly one interior lattice point.
			Then $P$ admits an affine dimer.
		\end{proposition}
		\begin{proof}
			There are 16 equivalence classes of convex lattice polygons with exactly one interior point \cite{poonenLattice12}, \cite{genus0polygons}.
			As seen in Figure \ref{figure:genus1easyClasses}, four of them are triangles, which admit an affine dimer by Proposition \ref{prop:triangles}.
			Three of them consist of pairwise antiparallel side segments, so admit an affine dimer by Proposition \ref{prop:doubleEverything}.
			Another three of them are obtained by adding a pair of antiparallel side segments $(\pm1, 0)$ to convex lattice polygons containing $(1,0)$ as a side segment which are already known to be dimers since they have no interior points. These admit an affine dimer by Proposition \ref{prop:addParallelEdges}.
			
			There are five equivalence classes left whose affine dimers are displayed in Figure \ref{figure:genus1checkedCases} \& \ref{figure:genus1checkedCases2}.
		\end{proof}
		
		\begin{figure}[H]\centering
			\includegraphics[width=4.5cm]{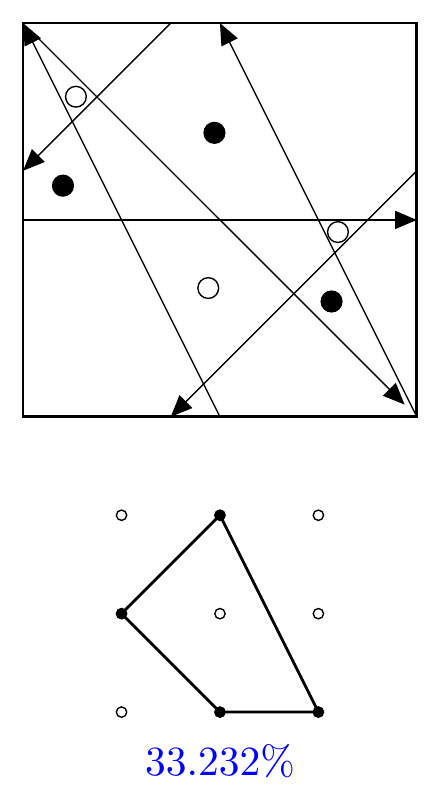}
			\hspace{1cm}
			\includegraphics[width=4.5cm]{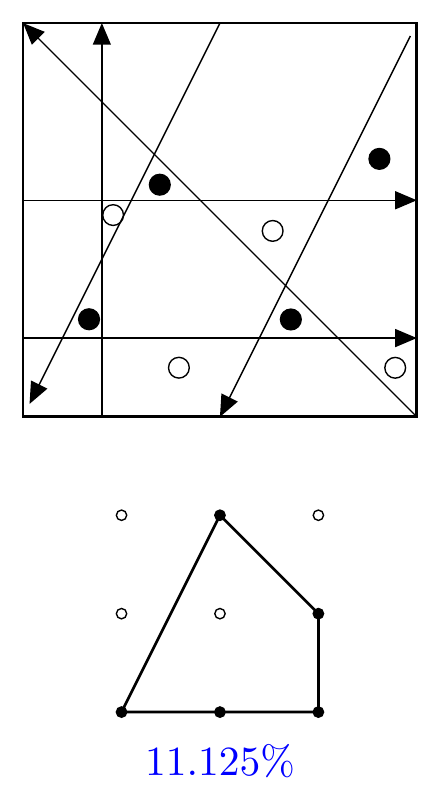}
			\hspace{1cm}
			\includegraphics[width=4.5cm]{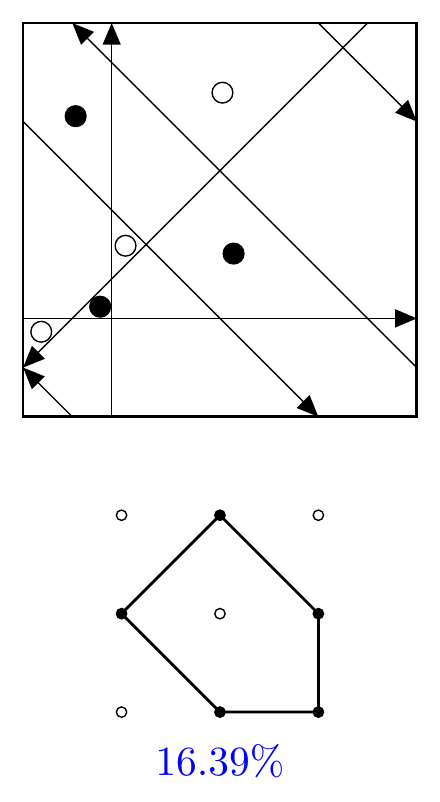}
			\caption{The five remaining equivalence classes of convex lattice polygons with one interior lattice point and affine dimers for them (continued in Figure \ref{figure:genus1checkedCases2}).
				The blue numbers are estimates ($\ge2\cdot 10^4$ trials) of the admissible volume.}
			\label{figure:genus1checkedCases}
		\end{figure}
		
		\begin{figure}[H]\centering
			\includegraphics[width=4.5cm]{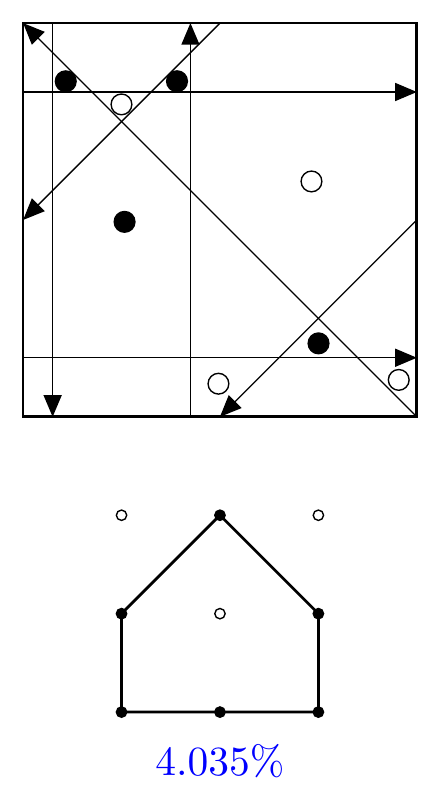}
			\hspace{1.4cm}
			\includegraphics[width=4.5cm]{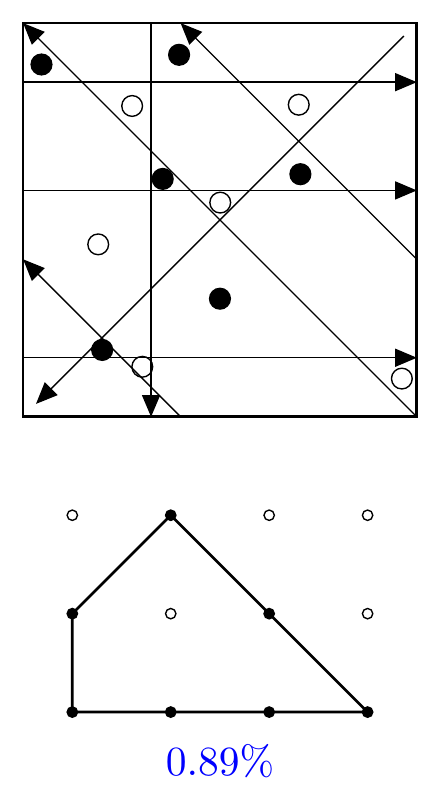}
			\caption{Continuation of Figure \ref{figure:genus1checkedCases}.}
			\label{figure:genus1checkedCases2}
		\end{figure}
		
		\subsection{Genus Two}
		
		The same holds for two interior points.
		
		\begin{proposition}
			\label{prop:genus2}
			Let $P$ be a convex lattice polygon with exactly two interior lattice points.
			Then $P$ admits an affine dimer.
		\end{proposition}
		\begin{figure}[h]\centering
			\includegraphics[width=5.5cm]{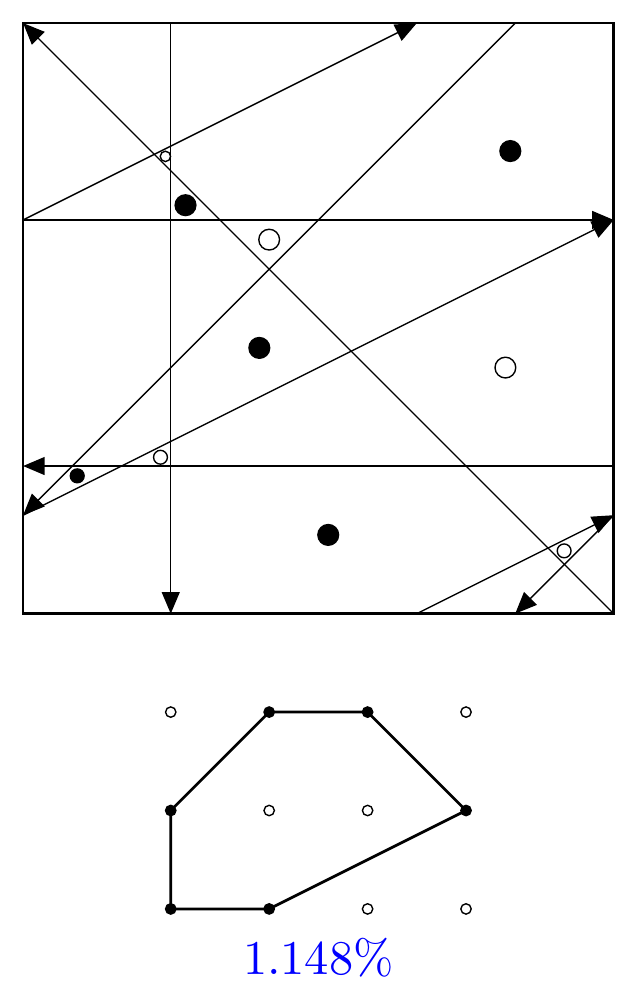}
			\caption{The only genus 2 hexagon whose dimer cannot be constructed from a lower genus dimer using the constructions of Section \ref{sec:new_dimers_from_old}.
				The blue number is an estimate ($10^5$ trials) of the admissible volume.}
			\label{figure:genus2example}
		\end{figure}
		\begin{proof}
			A classification up to equivalence of convex lattice polygons with two interior lattice points is provided by \cite{genus2polygons}.
			
			There are five classes of triangles, all of which admit an affine dimer by Proposition \ref{prop:triangles}.
			
			There are 19 classes of quadrilaterals.
			Three of them are parallelograms and thus admit an affine dimer by Proposition \ref{prop:addParallelEdges}.
			Six of them are obtained by adding a pair of antiparallel edges parallel to an existing edge to a convex lattice polygon of genus 0 or 1, and therefore admit an affine dimer by Propositions \ref{prop:genus0}, \ref{prop:genus1}, and \ref{prop:addParallelEdges}.
			The other 10 classes of quadrilaterals were checked manually to admit an affine dimer (see below).
			
			Similarly, all sixteen classes of pentagons and five classes of hexagons admit an affine dimer (see Figure \ref{figure:genus2example} for an example). There are no convex lattice $n$-gons with two interior lattice points and $n>6$.
			The 19 classes that required computer-aided verification can be found online at
			\url{https://jeffhicks.net/files/DHolmesSupplemental.pdf}.
			\todo{(B7) Add reference for the remaining genus 2 classes.}
		\end{proof}

		This completes the proof of Theorem \ref{thm:mainResult}.

		\section*{Acknowledgments} 
		
		The author would like to thank the \textit{London Mathematical Society} and the \textit{Department of Pure Mathematics and Mathematical Statistics, University of Cambridge} for their financial support in the form of an Undergraduate Research Bursary, and Jeff Hicks for suggesting and mentoring this project, and for providing invaluable feedback and many fruitful discussions.
		Finally, we would like to thank the referees for their helpful feedback, contributing multiple improvements to this paper.

		{\footnotesize  
			\medskip
			\medskip
			\vspace*{1mm} 
			
			\noindent {\it Daniel Holmes}\\  
			DPMMS, University of Cambridge\\
			Wilberforce Road\\
			Cambridge, CB3 0WB, United Kingdom\\
			E-mail: {\tt dh604@cam.ac.uk}\\ \\  
			
			
		
	}
	
	\vspace*{1mm}\noindent\footnotesize{\date{ {\bf Received}: October 7, 2021\;\;\;{\bf Accepted}: January 5, 2022}}\\
	\vspace*{1mm}\noindent\footnotesize{\date{  {\bf Communicated by Matthias Beck}}}
	
\end{document}